\catcode`\@=11

\magnification=1200
\hsize=121mm    \vsize=180mm
\hoffset=10mm    \voffset=10mm
\pretolerance=500  \tolerance=1000 
\brokenpenalty=500
\font\mathsy=msym10  
\font\rmgXVIII=cmbx9 at 18pt 
\font\rmgXV=cmbx9 at 15pt
\font\rmgXIII=cmbx9 at 13pt

\def\diagram#1{\def\normalbaselines
{\baselineskip=0pt          
\lineskip=10pt\lineskiplimit=1pt} 
\matrix{#1}}                

\def\hfl#1#2{\smash{\mathop{\hbox to 
10mm{\rightarrowfill}}                       
\limits^{\scriptstyle#1}_{\scriptstyle#2}}}                     

\def\gfl#1#2{\smash{\mathop{\hbox to
10mm{\leftarrowfill}}     
\limits^{\scriptstyle#1}_{\scriptstyle#2}}}               

\def\vfl#1#2{\llap{$\scriptstyle #1$}\left\downarrow
\vbox to 5mm{}\right. \rlap{$\scriptstyle #2$}}

\def\cro#1#2#3{\noindent{\bf [#1]} {\bf #2} : {\sl #3}}
{\nopagenumbers
\line{Vincent Mercat\hfil}\par
\line{138, rue Cardinet\hfil}\par
\line{75017 Paris\hfil}\par
\line{mercat@math.jussieu.fr\hfil}
\null\vskip 4cm
\centerline{\rmgXVIII Le probl\`eme }\bigskip
\centerline{\rmgXVIII de}\bigskip
\centerline{\rmgXVIII Brill-Noether}\bigskip
\centerline{\rmgXVIII pour des fibr\'es 
stables}\bigskip
\centerline{\rmgXVIII de petite pente.}\vskip 3cm

\vfill\break}

\noindent{\rmgXV Introduction}\vskip 0,8cm

\noindent Tout au long de ce travail, nous consid\'erons
que $C$ est une courbe alg\'ebrique lisse de genre $g$
sur un corps alg\'ebriquement clos de caract\'eristique
$0$.\medskip
\noindent Le probl\`eme de Brill-Noether consistait  \`a
l'origine \`a  d\'eterminer \`a quelles conditions il existe un
fibr\'e vectoriel en droites $L$ sur $C$, de
degr\'e
$d$ et poss\'edant $k$ sections globales ind\'ependantes.
On obtenait ainsi un  morphisme de la courbe $C$ dans
l'espace projectif {\mathsy P}$^{k-1}$ et la r\'esolution
du probl\`eme de Brill-Noether dans ce cas  a permis de
faire consid\'erablement avancer la classification des
courbes et l'\'etude de la vari\'et\'e de modules de
courbes ${\cal M}_g$. Les premiers r\'esultats remontent
\`a la fin du si\`ecle dernier. On trouve une approche
moderne de ce probl\`eme avec une d\'emonstration de
tous les th\'eor\`emes connus  dans l'ouvrage de
Arbarello, Cornalba, Griffiths et Harris (cf [{\bf
A-C-G-H}]). \par
\noindent Avec l'apparition des vari\'et\'es de modules
$U_{n,d}$ (resp. ${\widetilde U}_{n,d}$) de fibr\'es
vectoriels stables (resp. semi-stables) de rang $n$ et de
degr\'e $d$, le probl\`eme de Brill-Noether se
g\'en\'eralisait: \`a quelles conditions existe-il des fibr\'es
stables de rang $n$ de degr\'e $d$ poss\'edant $k$
sections globales ind\'ependantes?\medskip
\noindent En g\'en\'eralisant les constructions faites
dans le cas des fibr\'es en droites, on est amen\'e \`a
poser (voir le {\bf chapitre 1}):
$$W_{n,d}^{k-1}:=\{E\in U_{n,d}\mid h^0(E)\geq
k\}$$
$${\widetilde W}_{n,d}^{k-1}:=
\{[E]\in\widetilde{U}_{n,d}\mid h^0(\hbox{gr}\,E)\geq
k\}\ .$$
L'apparition de l'exposant $k-1$ est due au cas des
fibr\'es en droites: on cherchait un morphisme
$C\rightarrow \hbox{\mathsy P}^{k-1}$. Mais pour les
fibr\'es de rang sup\'erieur, cette notation est
particuli\`erement malheureuse.\par
\noindent On donne 
\`a
$W_{n,d}^{k-1}$ (resp. ${\widetilde W}_{n,d}^{k-1}$) une
structure de sous-sch\'ema ferm\'e de $U_{n,d}$ (resp. 
$\widetilde{U}_{n,d}$)
(cf [{\bf L}], [{\bf A-C-G-H}] ou [{\bf B-G-N}]).
Le probl\`eme de Brill-Noether consiste alors \`a trouver
quand ces espaces sont vides, quelle est leur
dimension, s'ils sont irr\'eductibles etc...\par
\noindent Il r\'esulte de la construction que si
$W_{n,d}^{k-1}\not = \emptyset$ et si
$W_{n,d}^{k-1}\not = U_{n,d}$, alors
 on sait que
$$\hbox{dim}\bigl(W_{n,d}^{k-1}\bigr)\geq
\rho(g,d,n,k-1)= n^2(g-1)+1-k(k-d+n(g-1))$$
Jusqu'au d\'ebut des ann\'ees 90, nous avions tr\`es
peu de r\'esultats (cf [{\bf Se}]) pour les fibr\'es de rang
sup\'erieur \`a un. Monserrat Teixidor i Bigas d\'emontre
en 1991 un th\'eor\`eme qui donne une solution assez
g\'en\'erale au probl\`eme (cf [{\bf Te1}] ou th\'eor\`eme
B-2 du Chap. 1). Ce th\'eor\`eme est encore aujourd'hui
le r\'esultat le plus g\'en\'eral que nous poss\'edons.
Cependant il n'est pas enti\`erement satisfaisant: il n'est
valable que pour une courbe g\'en\'erique (dans un sens
tr\`es vague) et des r\'esultats plus r\'ecents montrent 
qu'il ne donne pas une solution optimale. En 95,
Brambila-Paz, Grzegorczyk et Newstead (cf [{\bf B-G-N}] ou
th\'eor\`eme B-4 du chap. 1) donne une solution au
probl\`eme pour toutes les courbes si ${d\over
n}\leq1$.\medskip \noindent Ce travail
consiste en fait \`a g\'en\'eraliser ce th\'eor\`eme au cas des
fibr\'es stables de pente $1<\mu<2$.   Nous nous inspirons
tr\`es largement des id\'ees contenues dans [{\bf
B-G-N}] pour d\'emontrer le th\'eor\`eme suivant:\bigskip
\noindent{\bf Th\'eor\`eme:} {\it Soit $C$ une courbe
lisse. On suppose que $1<{d\over n}<2$. Alors les
espaces de Brill-Noether ${ W}_{n,d}^{k-1}$ sont non-vides
si et seulement si $$k\leq n+{1\over
g}(d-n)$$ et alors toutes les composantes irr\'eductibles de ${
W}_{n,d}^{k-1}$ sont de dimension exactement
$\rho(g,d,n,k-1)$.}\medskip
\noindent Pour
d\'emontrer ce th\'eor\`eme, nous commen\c cons dans
le {\bf chapitre 2}  par
\'etudier
les fibr\'es vectoriels engendr\'es par leurs
sections: si $F$ est un fibr\'e vectoriel de rang $l$ et
de degr\'e  $d$ engendr\'e par ses sections, alors on
note $D(F)$ le dual du noyau du morphisme
d'\'evaluation:
$$O\hfl{}{}D(F)^*\hfl{}{} {\rm H}^0(F)\otimes {\cal O}\hfl{}{}
F\hfl{}{}0$$
Nous obtenons que si $F$ est stable de rang $l$
et de pente ${d\over l}>2g$, $D(F)$ est un fibr\'e
stable de rang $n=d-gl$ et de degr\'e $d$ poss\'edant
$n+l$ sections globales ind\'ependantes. La condition
${d\over l}>2g$ implique que ${d\over n}<2$. On
obtient  un isomorphisme \par
\centerline{$W_{n,d}^{n+l-1}\simeq U_{l,d}$,}
\noindent ceci donne un plongement\par
\centerline{$U_{l,d}\hookrightarrow U_{d-gl,d}$.}
\noindent A notre connaissance, l'existence  de tels
morphismes entre des vari\'et\'es de modules de fibr\'es
stables \'etaient inconnus jusque l\`a: nous y
consacrerons ult\'erieurement une \'etude.\par
\noindent De plus, Paranjape et
Ramanan (cf [{\bf P-R}]) ont montr\'e que si $K$ est
le fibr\'e canonique, alors $D(K)$, que nous noterons
$E_K$, est un fibr\'e stable si la courbe $C$ est
non-hyperelliptique (semi-stable sinon). On en d\'eduit la
non-existence de fibr\'e stables de pente
$<2$ avec $k>n+{1\over g}(d-n)$. Nous
pr\'eciserons quelques propri\'et\'es
remarquables de ce fibr\'e.\medskip

\noindent {\bf Le chapitre 3} est le chapitre technique de
ce travail: le chapitre 2 donne l'existence de fibr\'es
stables tels que le rang et le degr\'e v\'erifient $d-n=gl$;
il reste donc \`a traiter le cas des fibr\'es de rang $n$ et
de degr\'e $d$ tels que $d-n$ n'est pas divisible par $g$.
Nous sommes conduits \`a faire
des calculs de param\`etres un peu lourds et
complexes, mais nous n'avons pas trouv\'e
d'\'echappatoire.\medskip
\noindent Notons que ce th\'eor\`eme montre en fait que
si ${d\over n}\geq1$, ${k\over
n}\leq1$ et $({d\over n },{k\over n})\not = (a,1)$ avec
$a$ un entier, alors
$W_{n,d}^{k-1}$ est non vide. \medskip

\noindent Je remercie J.M. Drezet, I. Grzegorczik, Y. Laszlo,
P. E. Newstead et C. Sorger pour m'avoir soutenu dans mon
travail.\vskip 1cm

\noindent{\rmgXV Chapitre 1: Le probl\`eme de
Brill-Noether}\vskip 0,8cm

\noindent Dans ce chapitre nous d\'efinissons les
espaces de Brill-Noether (partie A), puis  nous
donnons un panorama des diff\'erents r\'esultats
connus jusqu'\`a aujourd'hui (partie B). \vskip 0,7cm

\noindent{\rmgXIII A- Les espaces de Brill-Noether 
}\vskip 0,6cm 

\noindent Soit $C$ une courbe alg\'ebrique lisse 
de genre
$g\geq2$ sur un corps
alg\'ebriquement clos de caract\'eristique
$0$.\bigskip

\noindent{\bf Notations:}
\noindent On notera ${\cal O}={\cal O}_C$ le fibr\'e
trivial, $K$ le fibr\'e canonique sur $C$, et
$L^k=L\oplus\cdots \oplus  L$ (resp.
$L^{\otimes k}=L\otimes\cdots
\otimes  L$)  la somme (resp. le produit
tensoriel) de $k$ fois le fibr\'e $L$.\par

\noindent Un $g_d^k$ est un fibr\'e vectoriel 
en droites de
degr\'e $d$ poss\'edant $k+1$ sections globales
ind\'ependantes.\par
\noindent Si $E$ est un fibr\'e vectoriel 
alg\'ebrique de rang
$n$ et de degr\'e $d$ sur $C$, $\mu(E)={d\over n}$
est la {\it pente} de $E$. \par
\noindent Si $F$ est un fibr\'e vectoriel sur $C$, on
notera $n_F$ et $d_F$ le rang et le degr\'e de $F$. De
plus, on notera $D(F)$ le dual du noyau du morphisme
d'\'evaluation de $F$:
$$D(F)=Ker\bigl(ev : {\rm H}^0(F)\otimes {\cal
O}\longrightarrow F\bigr)^*\ .$$
\noindent On pose $h^i(E)=
\hbox{dim}\,{\rm H}^i(E)$.\medskip
\noindent On rappelle qu'un fibr\'e $E$  sur $C$ est
dit {\it stable} (resp. {\it semi-stable})  si pour tout
sous-faisceau propre
$F$ de
$E$ on a $\mu(F)<  \mu(E)$ (resp.
$\mu(F)\leq  \mu(E)$).\medskip
\noindent Tout fibr\'e $E$  semi-stable admet une
filtration
$$0=E_0\subset\cdots
\subset E_p=E$$
telle que $E_j/ E_{j-1}$ est stable 
et $\mu(E_j/ E_{j-1} )=\mu(E)$ pour $0<j\leq p $. Le
fibr\'e gradu\'e associ\'e $\oplus_j  E_j/ E_{j-1}$ est
appel\'e le {\it gradu\'e} de $E$ et est not\'e gr$\,E$.
La filtration sera appel\'ee {\it filtration de
Harder-Narasimhan}. Pour les fibr\'es non
semi-stables on a une filtration analogue
o\`u les quotients sont suppos\'es
semi-stables de pentes strictement
d\'ecroissantes. Cette filtration sera aussi
appel\'ee {\it filtration de
Harder-Narasimhan}.\medskip

\noindent On note $U_{n,d}$ (resp.
${\widetilde  U}_{n,d}$)  la vari\'et\'e de modules
des
fibr\'es vectoriels stables
(resp. semi-stables) de rang $n$
et de degr\'e $d$ sur $C$.  
Les points de $U_{n,d}$ sont les classes 
d'isomorphisme de fibr\'es stables tandis que les
points de ${\widetilde  U}_{n,d}$
repr\'esentent les classes
d'\'equivalences de fibr\'es 
semi-stables, o\`u $E\sim F$ si
et seulement si leurs gradu\'es 
gr$\,E$ et gr$\,F$ sont
isomorphes
(cf {\bf [Se]} ou{\bf [L]}) . On notera $[E]$ la classe
d'\'equivalence contenant $E$. \bigskip

\noindent{\bf Les espaces de Brill-Noether:} Les
espaces de Brill-Noether
$W_{n,d}^{k-1}$ et
${\widetilde W}_{n,d}^{k-1}$ sont d\'efinis
en tant qu'ensembles par:\par
$$W_{n,d}^{k-1}:=\{E\in U_{n,d}\mid h^0(E)\geq
k\}$$
$${\widetilde W}_{n,d}^{k-1}:=
\{[E]\in\widetilde{U}_{n,d}\mid h^0(\hbox{gr}\,E)\geq
k\}$$
On donne  \`a $W_{n,d}^{k-1}$
(resp. ${\widetilde W}_{n,d}^{k-1}$) une structure de
sous-sch\'ema ferm\'e de $U_{n,d}$ (resp. 
$\widetilde{U}_{n,d}$)
(cf [{\bf L}], [{\bf A-C-G-H}] ou [{\bf B-G-N}]).
Le probl\`eme de Brill-Noether consiste \`a trouver
quand ces espaces sont vides, quelle est leur
dimension, s'ils sont irr\'eductibles etc...\par
\noindent Il r\'esulte de la construction que si
$W_{n,d}^{k-1}\not = \emptyset$ et si
$W_{n,d}^{k-1}\not = U_{n,d}$, alors:\par 
- toute composante irr\'eductible $W$ de $W_{n,d}^{k-1}$
v\'erifie
$$\hbox{dim}\bigl(W\bigr)\geq
\rho(g,d,n,k-1)= n^2(g-1)+1-k(k-d+n(g-1))$$
\indent- son espace des points singuliers
Sing$\,W_{n,d}^{k-1}$ v\'erifie
$$W_{n,d}^{k}\subset\hbox{Sing}\,
W_{n,d}^{k-1}$$\par
- son espace tangent en un point $E$ 
tel que $h^0(E)=k$
est le noyau du morphisme
$$p^*:\hbox{Ext}^1(E,E)\hfl{}{}
{\rm H}^0(E)^*  \otimes{\rm
H}^1(E)$$ dual du morphisme de Petri
$$p:{\rm H}^0(E)
\otimes{\rm H}^0(E^*\otimes K)
\hfl{}{} {\rm H}^0(\hbox{End}(E)\otimes  K)$$
donn\'e par la multiplication des sections.\par
\noindent On  d\'eduit de cette derni\`ere propri\'et\'e
que si le morphisme de Petri est injectif 
en un tel point $E$, alors
$W_{n,d}^{k-1}$ est lisse en ce point et 
 que la composante
irr\'eductible contenant $E$ est de dimension
$\rho(g,d,n,k-1)$.\par
\noindent On ne sait pas en g\'en\'eral quand
 $W_{n,d}^{k-1}$ est non vide, irr\'eductible et quand
le
 morphisme de Petri est injectif. On sait
encore moins de choses sur ${\widetilde
W}_{n,d}^{k-1}$. En fait, ${\widetilde
W}_{n,d}^{k-1}$ n'est m\^eme pas la cl\^oture de 
$W_{n,d}^{k-1}$ dans ${\widetilde U}_{n,d}$ (cf  [{\bf
B-G-N}]).\medskip

\noindent Il est coutume de  poser $\mu={d\over
n}$ et $\lambda={k\over
n}$. L'id\'ee est due \`a A. King. Ceci permet de se
ramener
\`a un probl\`eme
\`a deux variables ($\mu$ et  $\lambda$) au lieu des
trois variables $d$, $n$ et $k$ et donne
une certaine analogie avec le cas des fibr\'es en
droites, comme nous le verrons par la suite.
De plus,  les r\'esultats vont pouvoir
s'interpr\'eter sur des graphiques ($\mu$ en
abscisse et $\lambda$ en ordonn\'ee). Si $E$ est
un fibr\'e sur $C$ de rang $n$ et de degr\'e $d$
poss\'edant $k$ sections globales ind\'ependantes,
on associe \`a $E$ le point de coordonn\'ees
$({d\over n},{k\over n})$, et le probl\`eme de
Brill-Noether devient: en un point de
coordonn\'ees $({d\over n},{k\over n})$ existe-t-il
des fibr\'es $E$ de rang $n$, de degr\'e $d$
poss\'edant $k$ sections globales
ind\'ependantes? Et quelle est la structure de
l'espace de Brill-Noether $W_{n,d}^{k-1}$
constitu\'es par tous ces fibr\'es?\medskip
\noindent Les th\'eor\`emes de 
Riemann-Roch et de Clifford permettent de
d\'egager une zone o\`u le probl\`eme n'est
pas trivial (cf figure a):\medskip

\noindent  {\bf Th\'eor\`eme de
Riemann-Roch:} {\it Soit $E$ un fibr\'e vectoriel de
degr\'e $d$ et de rang $n$ sur $C$, alors ${\rm
h}^0(E)-{\rm h}^1(E)=d+n(1-g)$. } D'o\`u ${\rm
h}^0(E)\geq d+n(1-g)$ ou encore $\lambda\geq
\mu+1-g$. Donc en dessous de la droite
$\lambda=\mu+1-g$, les espaces de
Brill-Noether sont les espaces $U_{n,d}$
tout entiers.\medskip

\noindent  {\bf Th\'eor\`eme de Clifford:} {\it
Soit
$E$ un fibr\'e semi-stable de rang $n$ et de degr\'e
$d$ sur $C$, tel que
$0\leq
\mu(E)\leq 2g-2$. Alors on a
$${\rm h}^0(E)\leq n+{d\over 2}\ .$$}
Or, $k\leq n+{d\over
2}\Leftrightarrow
\lambda\leq 1+{\mu\over 2} $.
Donc au-dessus de la droite d'\'equation $\lambda=
1+{\mu\over 2} $, les espaces de Brill-Noether sont
vides.\par
\noindent Notons que l'on trouve une d\'emonstration
du th\'eor\`eme de Clifford par G. Xiao dans [{\bf
B-G-N}] hormis la d\'emonstration des cas
d'\'egalit\'e que l'on peut lire dans  un article de 
R. Re (cf [{\bf Re}]).  Dans cet article l'auteur pr\'ecise
le th\'eor\`eme de Clifford. Nous donnons ici une
version simplifi\'ee des diff\'erents
th\'eor\`emes:\bigskip

\noindent  {\bf A-1 Th\'eor\`eme:} {\it Soit
$C$ une courbe non-hyperelliptique et soit $E$ un
fibr\'e vectoriel semi-stable de rang $n$ et de
degr\'e $d$ tel que
$1<\mu(E)<2g-2$. Alors on a
$$ {\rm h}^0(E)\leq {d+n\over
2}$$
qui peut aussi s'\'ecrire
$\lambda\leq{\mu+1\over 2}$.}\medskip

\noindent D'autre part, il est clair qu'un fibr\'e
semi-stable de degr\'e n\'egatif n'a pas de sections
globales.
Donc on peut se
restreindre
\`a 
$\mu>0$ et bien s\^ur \`a $\lambda\geq 0$.
De plus, si $E$ est un fibr\'e semi-stable
de pente $\mu>2g-2$, 
${\rm H}^1(E)={\rm H}^0(E^*\otimes K)=0$ puisque
la pente de $E^*\otimes K$ est strictement
n\'egative. D'apr\`es le th\'eor\`eme de
Riemann-Roch on obtient dans ce cas 
${\rm H}^0(E)=d-n(g-1)$. Donc si $\mu>2g-2$ tous
les fibr\'es se trouvent sur la droite de Riemann-Roch 
$\lambda=\mu-g+1$.\medskip

\noindent  On d\'eduit de ce qui pr\'ec\`ede une
premi\`ere zone d\'elimit\'ee par 
la droite de Riemann-Roch, celle de Clifford, les
axes et la droite $\mu =2g-2$ 
(cf figure a). La droite de Re n'intervient que
pour le cas non-hyperelliptique. En dehors
de cette zone, les espaces de Brill-Noether
sont soit vides soit la vari\'et\'e de modules
$U_{n,d}$ toute enti\`ere.\medskip

\noindent Enfin, pour que $W_{n,d}^{k-1}$ soit non
vide, on peut escompter que sa dimension
th\'eorique, $\rho(g,d,n,k-1)=
n^2(g-1)+1-k(k-d+n(g-1))$, soit positive. En fait cette
in\'egalit\'e ne s'exprime pas qu'en fonction des
variables $\mu$ et $\lambda$ et dans la pratique
l'in\'egalit\'e $\rho(g,d,n,k-1)\geq 1$ a un r\^ole
essentiel. D\'ej\`a elle s'exprime parfaitement
avec les seules variables 
$\mu$ et $\lambda$:\medskip
\centerline{$\rho(g,d,n,k-1)\geq1
\Leftrightarrow{1\over
n^2}\bigl(\rho(g,d,n,k-1)-1\bigr)\geq0$}
\noindent et\par 
\centerline{${1\over n^2}\bigl(\rho(g,d,n,k-1)-1\bigr)
=g-1-{k\over n}\bigl ({k\over n}-{d\over n}+g-1\bigr)
=g-\lambda(\lambda-\mu+g-1)-1$}\medskip

\noindent Donc\par
\centerline{$\rho(g,d,n,k-1)\geq1
\Leftrightarrow\rho(g,\mu,1,\lambda-1)
\geq1$}\medskip
\noindent On appelle la courbe d'\'equation
$\rho(g,\mu,1,\lambda-1) =1$, la {\it courbe de
Brill-Noether}. Cette \'equation est not\'ee dans de
nombreux articles ${\widetilde\rho}=0$.
L'\'equivalence ci-dessus montre que l'existence
d'un fibr\'e stable de pente ${d\over n} $ avec $k$
sections est li\'ee \`a l'existence de fibr\'es en droites
de "degr\'e" $\mu$ avec $\lambda$ sections (il
faudrait que ces valeurs soient enti\`eres): c'est dans
ce sens qu'il faut comprendre le th\'eor\`eme de
Teixidor (th\'eor\`eme B-2 de ce chapitre).\par
\noindent Sous cette courbe les espaces de
Brill-Noether sont escompt\'es
\^etre non vides (cf figure a).\par
\noindent Remarquons ici
que le th\'eor\`eme de Riemann-Roch et la dualit\'e
de Serre impliquent une certaine sym\'etrie
dans le graphique: si $E$ de pente $\mu$ a
$k$ sections, alors ${\rm h}^0(K\otimes E^*)
={\rm h}^1(E)=k+n(g-1)-d$. Donc si $E$
correspond \`a un point  $(\mu,\lambda)$,
alors $K\otimes E^*$ se place en $(2g-2-\mu,
\lambda+g-1-\mu)$. Il suffit donc de traiter le
cas $0\leq\mu\leq g-1$, le reste se
d\'eduisant par sym\'etrie.\medskip

\noindent   Le probl\`eme de
Brill-Noether se ram\`ene ainsi \`a l'\'etude de la
r\'egion que l'on a d\'egag\'ee pr\'ec\'edemment et
il faut d\'eterminer pour quels triplets $(n,d,k)$ de
cette r\'egion (et pour quelles courbes!) on a
bien:\medskip - $W_{n,d}^{k-1}$ (resp. ${\widetilde
W}_{n,d}^{k-1}$) est non vide\par
- $W_{n,d}^{k-1}$ (resp. ${\widetilde
W}_{n,d}^{k-1}$) est irr\'eductible, de  dimension
$\rho(g,d,n,k-1)$.\par
- Le lieu des points singuliers de $W_{n,d}^{k-1}$
v\'erifie
$${\rm Sing}\,W_{n,d}^{k-1} =W_{n,d}^{k}$$
\noindent Ceci a fait l'objet de nombreux
articles et l'on poss\`ede d\'ej\`a un certain nombre
de r\'esultats. Nous en donnons un panorama dans la
partie suivante.\vskip 0,7cm

\noindent{\rmgXIII B- Panorama des r\'esultats 
connus}\vskip 0,6cm

\noindent Pour les fibr\'es en droites, la
situation est particuli\`erement satisfaisante. Tous
les r\'esultats  se trouvent dans [{\bf A-C-G-H}]
chap.V et sont dus \`a une liste impressionnante
de math\'ematiciens. Les principaux r\'esultats dont 
nous aurons besoin par la suite  peuvent se mettre
sous la forme du   th\'eor\`eme suivant:\bigskip 

\noindent {\bf B-1 Th\'eor\`eme:} {\it Si
$C$ est une courbe alg\'ebrique lisse de 
genre $g$ et si
$\rho(g,d,1,k-1)\geq0$, alors 
$W_{1,d}^{k-1}\not = \emptyset$.
Si $C$ est une courbe alg\'ebrique
lisse g\'en\'erique de
genre $g$, alors
$$W_{1,d}^{k-1}\not = \emptyset \Leftrightarrow 
\rho(g,d,1,k-1)\geq0$$
et si $\rho(g,d,1,k-1)\geq1$, on
a dim$\,W_{1,d}^{k-1}=
\rho(g,d,1,k-1)$ ainsi que
\hbox{Sing$\,W_{1,d}^{k-1}=
\,W_{1,d}^{k}$.}}\medskip

\noindent Pour les fibr\'es de rang sup\'erieur,
le th\'eor\`eme le plus g\'en\'eral
que nous connaissons est d\^u \`a M. Teixidor i Bigas
 (cf {\bf [Te-1]} 1990). La formulation en \'etait fort
complexe et n'a
\'et\'e
"simplifi\'ee" que ult\'erieurement (Grzegorczyk, King,
Newstead...) avec l'introduction des param\`etres
$\lambda$ et $\mu$. Ici,  nous en donnons une
version plus courte, mais  compl\`ete gr\^ace \`a la
sym\'etrie de Riemann-Roch signal\'ee
\`a la fin de la partie A (cf figures b et c):\bigskip

\noindent {\bf B-2 Th\'eor\`eme:} {\it Soit $C$ une
courbe lisse g\'en\'erique.\par 
\noindent Si pour des
valeurs enti\`eres  $\lambda_0$ et $\mu_0$,
$\rho(g,\mu_0,1,\lambda_0-1)> 1$, alors pour
toutes valeurs  $\mu={d\over n}\geq\mu_0$ et
$\lambda={k\over n}\leq \lambda_0$,
$W_{n,d}^{k-1}$ (resp. ${\widetilde
W}_{n,d}^{k-1}$) est non vide et poss\`ede une
composante irr\'eductible de la bonne dimension.\par
\noindent Si 
$\rho(g,\mu_0,1,\lambda_0-1)= 1$,
 l'assertion est encore vraie,
sauf pour  $\mu=\mu_0$ o\`u l'on doit se restreindre
aux fibr\'es semi-stables.
}\medskip

\noindent La d\'emonstration de Teixidor utilise
des techniques de d\'eg\'en\'eration de courbes. Elle
construit une vari\'et\'e de modules
${\cal U}_{n,d}$ (resp. ${\widetilde {\cal U}}_{n,d}$ sa
compactification) de familles de fibr\'es vectoriels
stables (resp. semi-stables) de rang
$n$ et de degr\'e $d$ index\'ee par une famille de
courbes "locale": on consid\`ere une famille
projective et plate de courbes $\pi:  X\rightarrow S$ 
avec $X$ lisse et $S$ le spectre d'un anneau de
valuation discr\`ete (complet ou Hens\'elien cf [{\bf
E-H-1}] p. 347 et suivantes) telle que la fibre
g\'en\'erique soit non singuli\`ere et la fibre sp\'eciale
n'ait comme singularit\'e que des points doubles
ordinaires. Une famille de fibr\'es stables (resp.
semi-stables)  est alors un faisceau coh\'erent sur
$X$ tel que sa  restriction \`a la  fibre g\'en\'erique
soit un fibr\'e stable (resp. semi-stable) et que sa
restriction \`a la fibre sp\'eciale soit un faisceau de
profondeur 1 $b$-stable (resp. $b$-semi-stable) (cf 
[{\bf Se}] Chap. 7). Ceci n\'ecessite de polariser la
courbe sp\'eciale puisqu'elle n'est pas lisse. La
construction ne diff\`ere alors en rien de celle
donn\'ee dans [{\bf Se} ] Chap. 8. Puis elle suppose
que la fibre sp\'eciale
$T$ est compos\'ee d'une courbe  rationnelle lisse
intersect\'ee par $g$ courbes elliptiques. La partie
technique est alors de donner des conditions pour
qu'un fibr\'e $b$-stable sur
$T$ avec
$k$ sections se prolonge en une famille de fibr\'es
sans perdre de sections. C'est en fait la
g\'en\'eralisation des s\'eries lin\'eaires limites
introduites par Eisenbud et Harris (cf  [{\bf
E-H-1-2-3}]).\par
\noindent Notons que la vari\'et\'e de
modules ${\widetilde {\cal U}}_{n,d}$
\'etant projective, un fibr\'e semi-stable ${\cal F}_g$
sur la fibre g\'en\'erique se prolonge en un faisceau
${\cal F}$  sans torsion sur $X$ tel que sa
restriction ${\cal F}_0$ \`a la fibre sp\'eciale soit un
faisceau de profondeur 1. Aux points singuliers,
${\cal F}_0$ peut ne pas \^etre localement libre. Il
faut alors utiliser des blowing up pour obtenir un
fibr\'e au-dessus de la courbe sp\'eciale. Mais si l'on
suppose que la courbe sp\'eciale est lisse, ${\cal
F}_0$ est imm\'ediatement un fibr\'e et le
th\'eor\`eme de semi-continuit\'e implique que 
${\rm h}^0({\cal F}_0)\geq {\rm h}^0({\cal F}_g)$ (cf
[{\bf H}] th\'eor\`eme 12.8 p. 288). On en d\'eduit le
lemme suivant:\bigskip

\noindent {\bf B-3 Lemme:} {\it Soit $g,n,d,k$ des
entiers tels que, pour une courbe lisse g\'en\'erique
de genre
$g$, il existe un fibr\'e semi-stable $E$ de rang $n$,
de degr\'e $d$ et poss\'edant au moins $k$ sections
globales. Alors l'existence d'un tel $E$ est vraie pour
toutes les courbes lisses.}\medskip

\noindent {\it D\'emonstration:} 
 Il faut alors comprendre le terme
g\'en\'erique comme "il existe un ouvert $V$ de la
vari\'et\'e de modules de courbes ${\cal M}_{g}$ tel
que sur toute courbe lisse dans $ V$ il existe un bon
fibr\'e". Soit $C$ une courbe lisse, alors il existe une
famille projective et plate de courbes $\pi: 
X\rightarrow S$ comme ci-dessus telle que la fibre
g\'en\'erique g\'eom\'etrique soit dans $V$ et la fibre
sp\'eciale \'egale \`a $C$.  Alors le fibr\'e semi-stable
$E$ sur la fibre g\'en\'erique se prolonge \`a $X$ pour
donner un fibr\'e semi-stable sur la courbe $C$ et il
poss\`ede au moins autant de sections globales que
$E$.\hbox to 1cm{}$\diamondsuit$
\medskip

\noindent Ce lemme est faux pour les fibr\'es
stables: la vari\'et\'e ${ {\cal U}}_{n,d}$ n'est pas
compacte. Le lemme A-3 du chapitre 2 (cf plus loin)
montre par exemple que pour $C$
non-hyperelliptique il existe un fibr\'e stable, not\'e
$E_K$, de rang
$g-1$ de degr\'e $2g-2$ poss\'edant $g$ sections
globales alors que si la courbe est hyperelliptique, il
n'existe que des fibr\'es semi-stables  non stables
qui poss\`edent au moins $g$ sections. Ce cas
n'est pas
unique.\medskip
\noindent L'existence de fibr\'e stables dans le
th\'eor\`eme de Teixidor, comme on vient de le voir,
n'est valable que pour une courbe g\'en\'erique.
Par contre, pour les fibr\'es de petite
pente, de pente 
$\leq1$, on trouve une description  pr\'ecise et
valable pour toutes les courbes des espaces de
Brill-Noether de fibr\'es stables. Le th\'eor\`eme
suivant est d\'emontr\'e dans   [{\bf B-G-N}]
(cf graphique d):\bigskip

\noindent {\bf B-4 Th\'eor\`eme:} {\it Soit $C$ une
courbe lisse. On suppose que $0\leq{d\over n}\leq1
$.\par
\noindent $W_{n,d}^{k-1}$ est non vide si et
seulement si \par
\centerline{$d>0$, $n\leq d+(n-k)g$ et $(n,d,k)\not
=(n,n,n)$}\par \noindent et alors $W_{n,d}^{k-1}$ est
irr\'eductible
de la bonne dimension et v\'erifie\par
\centerline{Sing$\,W_{n,d}^{k-1}=W_{n,d}^{k}$ }\par
\noindent
\noindent ${\widetilde
W}_{n,d}^{k-1}$ est non vide si et
seulement si \par
\centerline{$\bigl(d=0$ et $k\leq n\bigr)$
 ou $\bigl(
d>0$ et $n\leq d+(n-k)g\bigr)$}\par
\noindent et alors ${\widetilde
W}_{n,d}^{k-1}$ est irr\'eductible.}\medskip

\noindent L'in\'egalit\'e $n\leq d+(n-k)g$
peut s'\'ecrire $1\leq\mu+ g(\mu-\lambda)$, ou
encore $\lambda\leq 1+{1\over g}(\mu-1)$. On
notera $\Delta$ la droite d'\'equation $\lambda=
1+{1\over g}(\mu-1)$. La droite $\Delta$ est en fait
la tangente \`a la courbe de Brill-Noether 
$\rho(g,\mu,\lambda,k-1) =1$ au point
$(1,1)$. Pour $0\leq \mu\leq 1$, les espaces
de Brill-Noether correspondants \`a des points
au-dessus de cette droite sont vides.\par
\noindent De plus, la clef de vo\^ute de la
d\'emonstration est: si
$E$ est un fibr\'e semi-stable de pente $\leq 1$ avec
$k$ sections globales ind\'ependantes, alors il peut
s'\'ecrire comme une extension de fibr\'es
$$0\hfl{}{} {\cal O}^k\hfl{}{} E\hfl{}{} F\hfl{}{}0\ .$$
Et on  v\'erifie  que ces extensions sont
param\'etr\'ees par une vari\'et\'e irr\'eductible de
dimension attendue (aux automorphismes pr\`es de
${\cal O}^k$). On voit que ce th\'eor\`eme donne
une solution compl\`ete pour les courbes de genre
$2$ grace \`a la sym\'etrie de Riemann-Roch (cf
figure e). C'est la raison pour laquelle dans les
chapitres suivants, nous supposerons que $g\geq
3$, m\^eme si les r\'esultats restent le plus souvent
valables dans le cas $g=2$.\par

\noindent En fait, le travail qui suit a pour
objet principal de montrer que l'in\'egalit\'e ci-dessus est
encore valable pour une pente $\mu$ comprise entre
$1$ et
$2$ (cf figure f).\medskip

\noindent Il existe de nombreux autres r\'esultats
moins g\'en\'eraux comme par exemple:\par
- Si $0<d\leq g-1$, $W_{n,d}^{0}$ est irr\'eductible
de la bonne dimension (cf [{\bf Su}]) et
Sing$\,W_{n,d}^{0}=W_{n,d}^{1}$ (cf [{\bf L}]).\par
- Le cas des fibr\'es de rang $2$ est en partie trait\'e
dans [{\bf Su}], [{\bf Te-2-3}] et [{\bf T}].\par
- Les vari\'et\'es $W_{3,1}^{k-1}$ et $W_{3,2}^{k-1}$
sont d\'ecrites dans [{\bf B-N}].\par
\noindent L'\'etude de ces cas particuliers fait
appara\^\i tre que $W_{n,d}^{k-1}$ peut avoir une
dimension plus grande que $\rho$, peut ne pas
\^etre  r\'eduit et
 avoir un lieu de singularit\'e distincts de
$W_{n,d}^{k}$, et ceci m\^eme pour
$C$ g\'en\'erique.\vskip 1cm

\noindent{\rmgXV Chapitre 2: La droite $\bf \Delta$ 
}\vskip 0,8cm 
\noindent Nous montrons dans ce chapitre que pour
une pente $\mu$, $1<\mu<2$, les vari\'et\'es de
Brill-Noether ${\widetilde
W}_{n,d}^{k-1}$ sont vides si le point de
coordonn\'ees $({d\over n},{k\over n})$ se trouve 
au dessus de la droite
$\Delta$. De plus, pour tout point ayant un sens sur
cette droite, nous donnons une description
compl\`ete des espaces de Brill-Noether
correspondants. Il r\'esulte de la partie B que les
points sur la droite
$\Delta$ correspondent \`a des fibr\'es de rang
$n$ de degr\'e $n+gl$ poss\'edant $n+l$ sections
globales, o\`u $l$ est un entier. Dans la partie C
nous donnons une description g\'en\'erale des espaces de
Brill-Noether dans le cas o\`u ceux-ci sont non-vides. Il
nous restera \`a prouver l'existence de fibr\'es stables
correspondant \`a des points sous la droite $\Delta$ quand
$d=n+gl+l'$ avec $0< l'<g$. Ceci sera fait dans le
chapitre 3.\medskip

\noindent La partie A traite la non-existence.\par
\noindent La partie B traite l'existence et d\'ecrit les
espaces de Brill-Noether correspondants pour
les points qui se trouvent sur la droite
$\Delta$.\par
\noindent La partie C pr\'ecise la structure des espaces de
de Brill-Noether pour ${d\over n}\leq2$.
\vskip 0,7cm

\noindent{\rmgXIII A- Non-existence de fibr\'es
stables.}\vskip 0,6cm

\noindent Rappelons tout d'abord que la droite
$\Delta$ a pour \'equation $\lambda=1+{1\over g}
(\mu-1)$. Le th\'eor\`eme ci-dessous implique qu'au-dessus
de cette droite les espaces de Brill-Noether sont
vides.\bigskip

\noindent{\bf A-1 Th\'eor\`eme:} {\it Si $E$ est un
fibr\'e semi-stable sur $C$ de rang $n$ et de pente
$\mu$,
$0<\mu={d\over n}<2$, alors
$${\rm h}^0(E)\leq n+{1\over g}(d-n)$$
ce qui s'\'ecrit aussi
$$\lambda\leq1+{1\over g}
(\mu-1)$$
o\`u $k={\rm h}^0(E)$ et $\lambda={k\over
n}$.}\medskip

\noindent Par la suite, nous aurons besoin
d'un r\'esultat plus g\'en\'eral. L'hypoth\`ese
de la semi-stabilit\'e du fibr\'e $E$ est
souvent difficile \`a v\'erifier. La proposition
ci-dessous  permet de contourner ce
probl\`eme dans certains cas. Le
th\'eor\`eme  A-1 s'en d\'eduit
imm\'ediatement.\par

\noindent {\bf A-2 Proposition:} {\it Soit $E$ un fibr\'e 
vectoriel de rang $n$ et de degr\'e $d$
tel que:\medskip
- son sous-fibr\'e semi-stable maximal est de pente
$<2$;\par
- h$^0(E^*)=0$.\medskip
\noindent Alors 
$${\rm h}^0(E)\leq n+{1\over g}(d-n)\ \ \ \ (**) \ \ .
$$}
\noindent Nous aurons besoin du lemme suivant qui
est d\^u \`a Paranjape et Ramanan (cf [{\bf
P-R}]):\bigskip
\noindent{\bf A-3 Lemme:} {\it Soit
$K$ le fibr\'e canonique, alors $D(K)$, le fibr\'e vectoriel de
rang $g-1$ et de degr\'e
$2g-2$ dual du noyau du morphisme d'\'evaluation ${\rm
H}^0(K)\otimes {\cal O}\rightarrow K$:
$$0\hfl{}{} D(K)^*\hfl{}{}
{\rm H}^0(K)\otimes {\cal
O}\hfl{}{}K\hfl{}{}0$$
 est stable si la courbe $C$ n'est pas
hyperelliptique, sinon on a \hbox{$D(K)\simeq
L\oplus
\cdots
\oplus L$} o\`u $L$ est l'unique $g_2^1$ de $C$.
Dans ce cas, $D(K)$ est 
semi-stable, non stable.}\medskip
\noindent Nous noterons $E_K$, le fibr\'e $D(K)$. Enfin, pour
les fibr\'es de pente $2$,  on a le corollaire suivant de la
proposition  A-2:\bigskip
\noindent{\bf A-4 Corollaire:} {\it
Si $E$ est stable de pente 2,  l'in\'egalit\'e $(**) $
est toujours valable pour $E
\not\simeq E_K$ et $C$ non hyperelliptique ou $E
\not\simeq L$ et $C$ hyperelliptique de 
$g_2^1$ isomorphe \`a $L$.}\medskip

\noindent{\it D\'emonstration de la proposition A-2 et
de son corollaire:} Consid\'erons la suite exacte
$$0\hfl{}{} E_K^*\hfl{}{}{\rm H}^0(K)\otimes {\cal
O}\hfl{}{}K\hfl{}{}0$$
On tensorise cette suite par un fibr\'e $E$
$$0\hfl{}{} E_K^*\otimes E\hfl{}{} {\rm H}^0(K)
\otimes E\hfl{}{} K\otimes
E\hfl{}{}0
$$
pour obtenir une suite exacte longue:
$$0\longrightarrow {\rm H}^0(E\otimes E_K^*)
\longrightarrow {\rm H}^0(E)\otimes{\rm H}^0(K)
\longrightarrow {\rm H}^0(E\otimes K)$$
$$\longrightarrow {\rm H}^1(E\otimes E_K^*)
\longrightarrow {\rm H}^1(E)\otimes{\rm H}^0(K)
\longrightarrow {\rm H}^1(E\otimes K)
\longrightarrow0$$
Supposons que ${\rm h}^0(E\otimes E_K^*)
=0$ et que
${\rm h}^1(E\otimes K)=0$.
 Alors d'une part, 
on a h$^0(E\otimes K)\geq gk$ avec h$^0(E)=k$
et d'autre
part, d'apr\`es le th\'eor\`eme de
Riemann-Roch, on a\medskip
\line{\hbox to 3cm{\hfil\hfil${\rm h}^0(E\otimes
K)$\hfil}$=
{\rm h}^1(E\otimes K)
+n(2g-2)+d-n(g-1)$\hfil}\par
\line{\hbox to
3cm{}$= ng-n+d$\hfil}\medskip
\noindent  Donc, 
$$kg\leq ng-n+d$$
d'o\`u
$$k\leq n+{d-n\over g}$$
Ce qui donne en divisant par $n$ et en posant
$\mu={d\over n}$ et $\lambda={k\over n}$:
$$\lambda\leq1+{1\over g}
(\mu-1)$$
Ce que nous voulions.\medskip
\noindent Il reste donc \`a montrer que dans les
conditions de la proposition on a bien
h$^0(E\otimes E_K^*)={\rm h}^1(E\otimes
K)=0$.\par
\noindent Si le sous-fibr\'e maximal de $E$ est de
pente $<2$, alors il n'existe pas de morphisme non
nul de $E_K$ semi-stable de pente $2$ dans $E$,
donc h$^0(E_K^* \otimes E)=0$. Si $E$ est stable de
pente $2$, le m\^eme raisonnement est valable dans
les conditions du corollaire.\par
\noindent D'apr\`es le th\'eor\`eme de dualit\'e de
Serre, on a: h$^1(E\otimes K)={\rm h}^0(E^*)$.
Par hypoth\`ese  h$^0(E^*)=0$, ce qui termine la
d\'emonstration.\hbox to
0,5cm{}$\diamondsuit$\medskip

\noindent{\bf A-5 Remarque:} Supposons
$d>0$. On a montr\'e en fait qu'en un point $({d\over
n},{k\over n})$
 au-dessus de la droite $\Delta$, il ne peut exister
de fibr\'e semi-stable $E$ de rang $n$, de
degr\'e $d$ et poss\'edant  $k$ sections
ind\'ependantes que si h$^0(E\otimes
E_K^*)\not=0$ et que h$^0(E\otimes E_K^* )$
mesure en fait le nombre  maximal de sections 
suppl\'ementaires possibles.\medskip
\noindent Pour majorer le nombre de sections
globales in\'ependantes de $E$, on voit qu'il sera
int\'eressant de majorer h$^0(E\otimes E_K^*)$. Dans
ce but, on donne  pour $E_K$ une sorte de lemme de
dualit\'e de Serre\bigskip
\noindent{\bf A-6 Lemme:} {\it Soit $E$ un fibr\'e
vectoriel, alors ${\rm H}^0(E\otimes E^*_K)\simeq{\rm
H}^1(D(E))^*$.}\medskip
\noindent{\it D\'emonstration:} On a vu que
H$^0(E\otimes E_K^*)$ est le noyau du morphisme ${\rm
H}^0(E)\otimes {\rm H}^0(K)\rightarrow {\rm
H}^0(E\otimes K)$ donn\'e par la multiplication des
sections. Or, la suite exacte 
$$0\hfl{}{}D(E)^*\otimes K\hfl{}{}{\rm H}^0(E)\otimes
K\hfl{ev\otimes Id}{} E\otimes K$$
montre que ce noyau est encore isomorphe \`a
H$^0(D(E)^*\otimes K)$, ce qui termine la
d\'emonstration du lemme.\hbox to 1cm{}
$\diamondsuit$
\vskip 0,7cm 

\noindent{\rmgXIII B- Les points de la droite $\bf
\Delta$}\vskip 0,6cm

\noindent Un fibr\'e vectoriel stable (resp.
semi-stable) $E$ de pente $\mu$, $1<\mu<2$
correspondant \`a un point sur la droite $\Delta$ doit
v\'erifier $$kg=d+n(g-1)$$
On posera avantageusement $k=n+l$.  Alors
$$kg=d+n(g-1)\Leftrightarrow g(n+l)=d+n(g-1)
\Leftrightarrow d=n+gl$$
De plus, $1<\mu(E)={d\over n}=1+{gl\over n} <2\,$
implique 
$$0<{gl\over n}<1\ ,$$
ce qui montre que $l>0$ et $n>gl$.\medskip
\noindent Les espaces de Brill-Noether que l'on se
propose de d\'ecrire ici sont donc
$W_{n,n+gl}^{n+l-1}$ (resp.
${\widetilde W}_{n,n+gl}^{n+l-1}$), pour $l>0$ et
$n>gl$ (cf figure f).\par
\noindent Rappelons que si $F$ un fibr\'e vectoriel, alors
on note
$$D(F)={\rm Ker}\bigl(ev:{\cal O}\otimes
{\rm H}^0(F)\hfl{}{}F\bigr)^*\ .
$$
\noindent Le th\'eor\`eme suivant donne la
description annonc\'ee:\bigskip

\noindent {\bf B-1 Th\'eor\`eme:} {\it Soit $n,l$ des
entiers, $l>0$ et $n>gl$. Alors on a des
isomorphismes
$$U_{l,n+gl}\hfl{\sim}{} W_{n,n+gl}^{n+l-1}$$
$$ F\mapsto D(F)$$
et
$${\widetilde U}_{l,n+gl}\hfl{\sim}{} 
{\widetilde W}_{n,n+gl}^{n+l-1}$$
$$ F\mapsto D(F)$$}\medskip

\noindent {\it D\'emonstration:} On donnera une
d\'emonstration du th\'eor\`eme essentiellement
pour les fibr\'es stables. Le cas semi-stable est
identique \`a quelques d\'etails pr\`es, que nous
nous efforcerons d'indiquer.\par

\noindent Soit $F$ un fibr\'e stable (resp.
semi-stable) de degr\'e $d=n+gl$ et de rang $l$. On
a ${d\over l}= {n+gl \over l}>2g$, puisque par
hypoth\`ese $n>gl$, $F$ est engendr\'e par ses
sections, h$^1(F)=0$ et h$^0(F)=d-l(g-1)= l+n$. On
pose $E\simeq D(F)$. On a une suite exacte
$$0\hfl{}{} E^*\hfl{}{} {\cal O}\otimes
{\rm H}^0(F)\hfl{}{} F\hfl{}{}0\ ,$$
et $E$
est un fibr\'e vectoriel de rang $n$, de degr\'e $d$ tel
que h$^0(E)\geq n+l$ et h$^0(E^*)=0$.
De plus, comme
$n>lg$,  $\mu(E)={d\over n}=1+{gl\over n}<2$.
L'\'enonc\'e du th\'eor\`eme est donc coh\'erent.
On va d'abord montrer que
$E$ est stable (resp. semi-stable), puis on
d\'efinit le morphisme
${ U}_{l,n+gl}\rightarrow 
{ W}_{n,n+gl}^{n+l-1}$ (resp. 
${\widetilde U}_{l,n+gl}\rightarrow 
{\widetilde W}_{n,n+gl}^{n+l-1}$) et pour terminer
on montre que c'est un isomorphisme (la
stabilit\'e de $E$ a d\'ej\'a \'et\'e d\'emontr\'ee
par
David Butler (cf [{\bf B}]).\medskip

\noindent {\bf Stabilit\'e de E}: On se restreint
maintenant au cas $F$ stable. 
Supposons que $E$ n'est
pas stable. Soit $G$ un fibr\'e quotient de $E$ stable de
plus petite pente, on a une suite exacte
$$0\hfl{}{}H\hfl{}{} E\hfl{}{} G\hfl{}{}0$$
o\` u $H$ est le noyau.
Par d\'efinition, on sait que
$\mu(G)\leq \mu(E)$ et donc $\mu(H)\geq
\mu(E)$.\medskip

\noindent 
D'apr\`es la proposition A-2, on sait que
$${\rm h}^0(G)\leq n_G+{1\over
g}(d_G-n_G)\ \ .$$
Le fibr\'e $G$ est engendr\'e par ses
sections, on peut poser h$^0(G)=n_G+l_G$,
o\`u $l_G$ est un entier strictement positif.
L'in\'egalit\'e pr\'ec\'edente s'\'ecrit
$$l_G\leq {1\over g}(d_G-n_G)$$
et donc
$$d_G\geq n_G+gl_G \ \ .$$
Or $d_G=d-d_H$ et $n_G=n-n_H$. L'in\'egalit\'e
devient
$$d-d_H\geq n-n_H+gl_G$$
et avec $d=n+gl$
$$g(l-l_G)+n_H\geq d_H$$
$\Leftrightarrow$\medskip
\line{\hfil ${g\over n_H}(l-l_G)+1\geq
{d_H\over n_H}$ .\hfil($*$)}\medskip
\noindent Par hypoth\`ese $\mu(H)\geq \mu(E)
=1+{gl\over n}>1$, donc $l>l_G$.\par 
\noindent D'autre part le diagramme
$$\matrix{
&&&&&&0&&&\cr
&&&&&&\downarrow&&&\cr
&& && &
& H&&\cr
&&&&&&\downarrow&&&\cr
0&\longrightarrow& F^*&\longrightarrow& {\cal
O}^{n+l}\simeq {\rm H}^0(F)^*\otimes {\cal O}&
\longrightarrow& E&\longrightarrow&0\cr
&&&&&&\downarrow&&&\cr
&& && &
& G&&\cr
&&&&&&\downarrow&&&\cr
&&&&&&0&&&\cr
}$$
peut se compl\'eter en un diagramme commutatif de
la mani\`ere suivante: on a un morphisme ${\rm
H}^0({\cal O}^{n+l})\rightarrow {\rm
H}^0(G)$ qui est la compos\'ee de ${\rm
H}^0({\cal O}^{n+l})\hookrightarrow {\rm
H}^0(E)$ et de ${\rm
H}^0(E)\rightarrow {\rm
H}^0(G)$. Soit $V$ l'image de ce morphisme;
comme ${\cal O}^{n+l}$ engendre $E$, on en
d\'eduit que $V$ engendre $G$ et donc que l'on
peut \'ecrire dim$\,V=n_G+l'_G$ avec $l'_G>0$. De
plus $V\hookrightarrow {\rm H}^0(G)$, donc $l'_G\leq
l_G<l$. Et le noyau du morphisme surjectif ${\cal
O}^{n+l} \rightarrow V\otimes {\cal O}$ est un fibr\'e
trivial de la forme ${\cal O}^{n_H+l'_H}$, c'est un
sous-espace de sections de ${\rm H}^0(H)$. Ce qui
donne le diagramme:

$$\matrix{
&&0&&0&&0&&&\cr
&&\downarrow&&\downarrow&&\downarrow&&&\cr
0&\longrightarrow&M^*&\longrightarrow& 
{\cal O}^{n_H+l'_H}&\longrightarrow
& H&&\cr
&&\downarrow&&\downarrow&&\downarrow&&&\cr
0&\longrightarrow& F^*&\longrightarrow& {\cal
O}^{n+l}& \longrightarrow&
E&\longrightarrow&0\cr
&&\downarrow&&\downarrow&&\downarrow&&&\cr
0&\longrightarrow&N^*&\longrightarrow& 
{\cal O}\otimes V&\longrightarrow
& G&\longrightarrow&0\cr
&&&&\downarrow&&\downarrow&&&\cr
&&&&0&&0&&&\cr
}$$
o\`u $M^*$
et $N^*$ sont les noyaux des morphismes
d'\'evaluation correspondant aux espaces de
sections. On a 
$l=l'_H+l'_G$, et comme $l_G\geq l'_G$,
on obtient $l'_H\geq l-l_G>0$. \par
\noindent De plus, $d_M\leq  d_H$: en effet,
$d_H-d_M $ est le degr\'e du coker de ${\cal
O}^{n_H+l'_H}\rightarrow H $ et $H$  a pour filtration de
Harder-Narasimhan $H_0\subset\cdots \subset
H_j=H$ et $\mu(H_j/ H_{j-1} )\geq\mu(G)$;  le
lemme ci-dessous permet de conlure (la
d\'emonstration en est triviale):\bigskip

\noindent {\bf B-2 Lemme:} {\it Soit $E$ un fibr\'e
vectoriel  de filtration de Harder-Narasimhan
$0=E_0\subset\cdots
\subset E_p=E$. Alors tout fibr\'e quotient de $E$
est de pente $\geq\mu(E_p/E_{p-1})$.}\medskip

\noindent  Retournons au diagramme: $F^*$ \'etant
stable, on doit alors avoir l'in\'egalit\'e
\hbox{$\mu(F)={d\over l}<\mu(M)={d_M\over
n_M}$} (notons que $M\not
\simeq0$ car $l'_H>0$ et $M\not \simeq F$ car
$N$ est non nul).  Or
$d_M\leq  d_H$ et $n_M\geq l'_H$, donc ${d\over
l}<\mu(M)\leq 
 {d_H\over l'_H}$. D'o\`u, avec $l'_H\geq (l-l_G)$
(cf plus haut),
 $$(l-l_G)\leq l'_H<d_H{l\over d}$$
Si l'on introduit cette in\'egalit\'e dans l'in\'egalit\'e
($*$), on obtient
$${gl\over d}{d_H\over n_H}+1>{d_H\over n_H}$$
$\Leftrightarrow$
$$1>(1-{gl\over d}){d_H\over n_H}$$
Or $d=n+gl$, donc $1-{gl\over d}={n\over d}$ et
l'in\'egalit\'e ci-dessus devient
$${d\over n}>{d_H\over n_H}$$
Ceci contredit l'hypoth\`ese $\mu(H)\geq\mu(E)$.
Donc $E$ est stable.\bigskip
\noindent Si $F$ est semi-stable le seul changement qui
intervient remplace la derni\`ere in\'egalit\'e par:
$${d\over n}\geq{d_H\over n_H}$$
et $E$ est alors
semi-stable.\bigskip

\noindent{\bf D\'efinition 
du morphisme $\bf { U}_{l,n+gl}\rightarrow 
{ W}_{n,n+gl}^{n+l-1}$:} On consid\`ere la vari\'et\'e
$R^s$ qui param\`etre une famille ${\cal E}$ de
fibr\'es de rang $l$ et de degr\'e $n+gl$, et qui sert
\`a construire la vari\'et\'e de modules $U_{l,n+gl}$
(cf [{\bf LP}] ou [{\bf Se}]):
le morphisme d\'eduit de ${\cal E}$
$$R^s\longrightarrow U_{l,n+gl}$$
est un bon quotient par un groupe du type
PGL($N$). Soit $\pi:R^s\times C\rightarrow R^s$ la
projection. Alors $\pi_*{\cal E}$ est un fibr\'e
vectoriel de rang $n+l$, et le morphisme
canonique sur $R^s\times C$
$$\pi^*\pi_*{\cal E}\hfl{\psi}{} {\cal E}$$
est surjectif. On pose
$${\cal F}={\rm Ker} (\psi)^*\ .$$
C'est une famille de fibr\'es stables de rang $n$ et
de degr\'e $n+gl$, param\'etr\'ee par $R^s$. On en
d\'eduit un morphisme  
$$R^s\hfl{}{} U_{n,n+gl}$$
qui est en fait dans $W_{n,n+gl}^{n+l-1}$. Ce
morphisme \'etant PGL($N$) invariant, on en
d\'eduit un morphisme (cf [{\bf M-F}]):
$$U_{l,d}\rightarrow W_{n,d}^{n+l-1}$$
Le m\^eme raisonnement est valable pour les
fibr\'es semi-stables et un
raisonnement analogue permettra de d\'efinir le
morphisme inverse.
\bigskip

\noindent{\bf Le morphisme $\bf {
U}_{l,n+gl}\rightarrow  { W}_{n,n+gl}^{n+l-1}$ est un
isomorphisme:} Pour  $n>gl$, on a donc par
d\'efinition: \par $$U_{l,n+gl}\hfl{}{}
W_{n,n+gl}^{n+l-1}$$
$$F\mapsto
D(F)$$
Ce morphisme est injectif car si l'on a deux suites
exactes
$$0\rightarrow E^* \rightarrow {\cal O}\otimes{\rm
H}^0(F) 
\rightarrow F \rightarrow 0$$
et
$$0\rightarrow E^* \rightarrow {\cal O}\otimes{\rm
H}^0(F')  
\rightarrow F' \rightarrow 0$$
alors h$^0(E)=n+l$ et le dual du noyau du
morphisme d'\'evaluation
$${\cal O}\otimes{\rm
H}^0(E) 
\rightarrow E \rightarrow 0$$
est isomorphe \`a $F$ et \`a $F'$, donc 
$F\simeq F'$.\medskip
\noindent Pour montrer que ce morphisme est un
isomorphisme, il faut montrer que, pour tout
$ E$ dans $ W_{n,n+gl}^{n+l-1}$,\medskip
-  $E$ est engendr\'e
par ses sections,\par
- Le noyau $F^*$ du morphisme d'\'evaluation 
${\cal O}^{n+l}\rightarrow E$ est stable.\medskip
\noindent En utilisant les propri\'et\'es universelles
des diff\'erents sch\'emas, le morphisme
r\'eciproque est alors donn\'e par
$$ E\mapsto D(E)$$
et il en est de m\^eme  pour les fibr\'es
semi-stables.\medskip

\noindent Soit donc $E\in W_{n,n+gl}^{n+l-1}$ et
supposons que $E$ n'est pas engendr\'e par ses
sections. Soit $Im$ le faisceau image du morphisme
d'\'evaluation
$${\cal O}\otimes{\rm H}^0(E) \hfl{}{} E$$
On pose ${Im}\simeq Im'\oplus {\cal O}^s$ o\`u
h$^0(Im'^*)=0$ (comme $Im$ est engendr\'e par ses
sections, c'est possible). On a
h$^0(Im')\geq n_{Im'}+l$.\par
\noindent De plus, tout sous-fibr\'e de $Im'$ est un
sous-faisceau de $E$ et donc de pente
$<\mu(E)<2$. On est donc dans les conditions de 
la proposition A-2, il en r\'esulte que:
$${\rm h}^0(Im')\leq n_{Im'}+{1\over
g}(d_{{ Im}'}-n_{Im'})$$
Comme ${\rm h}^0(Im')\geq n_{Im'}+l$, on obtient
$$n_{Im'}+l\leq n_{Im'}+{1\over
g}(d_{Im'}-n_{Im'})$$
dont on d\'eduit sans difficult\'e
$$\mu(Im')\geq  1+{gl\over n_{Im'}}\geq 
1+{gl\over n_E}=\mu(E)$$
et ceci contredit la stabilit\'e de $E$.
Donc $E$ est engendr\'e par ses sections. 
Pour $E$ semi-stable, on d\'eduit $n_{Im'}=n$, donc
$d_{Im'}=d$, d'o\`u $Im'\simeq E$. \medskip

\noindent Il reste \`a montrer que le fibr\'e $F=D(E)$
tel que
$$0\rightarrow F^* \rightarrow {\cal O}\otimes {\rm
H}^0(E) 
\rightarrow E \rightarrow 0$$
est un fibr\'e stable.\par
\noindent On montre tout d'abord que h$^0(F)=n+l$:
d'apr\`es le lemme A-6 on sait que h$^1(F)={\rm
h}^0(E\otimes E_K^*)=0$ et donc h$^0(F)=d-l(g-1)=n+l$.

\noindent Supposons maintenant que
$F$ n'est pas stable. On a alors une
suite exacte
$$0\hfl{}{}N\hfl{}{}F\hfl{}{} M\hfl{}{}0$$
telle que $N$ est semi-stable et 
${d_N\over n_N}=\mu(N)\geq\mu(F)={d\over l}>2g$.
$N$ est engendr\'e par ses sections et
h$^1(N)=0$. On en d\'eduit que
$${\rm H}^0(F)\simeq {\rm H}^0(N)\oplus
{\rm H}^0(M)$$
et donc un diagramme commutatif
$$\matrix{
&&0&&0&&0&&&\cr
&&\downarrow&&\downarrow&&\downarrow&&&\cr
0&\longrightarrow&G^*&\longrightarrow& 
{\cal O}\otimes{\rm H}^0(N)&\longrightarrow
& N&\longrightarrow&0\cr
&&\downarrow&&\downarrow&&\downarrow&&&\cr
0&\longrightarrow& E^*&\longrightarrow& {\cal
O}\otimes {\rm H}^0(F)& \longrightarrow&
F&\longrightarrow&0\cr
&&\downarrow&&\downarrow&&\downarrow&&&\cr
0&\longrightarrow&H^*&\longrightarrow& 
{\cal O}\otimes{\rm H}^0(M)&\longrightarrow
& M&\longrightarrow&0\cr
&&\downarrow&&\downarrow&&\downarrow&&&\cr
&&0&&0&&0&&&\cr
}$$
avec $H^*$ et $G^*$ les noyaux des morphismes
d'\'evaluation. Le rang de $G$, $n_G$, v\'erifie
$n_G={\rm h}^0(N)-n_N$. Mais on a vu que
h$^1(N)=0$, donc h$^0(N)=d_N-n_N(g-1)$. Il en
d\'ecoule que  $n_G=d_N-gn_N$.\par
\noindent Comme $E$ est stable, on doit avoir
$\mu(G)>\mu(E)$, ce qui s'\'ecrit
$${d_N\over d_N-gn_N}>{d\over n}$$
Or $n=d-lg$, l'in\'egalit\'e pr\'ec\'edente est
\'equivalente \`a
$$1+{gn_N\over d_N-n_Ng}>1+{gl\over d-lg}
\Leftrightarrow
{n_N\over d_N-n_Ng}>{l\over d-lg}$$
$\Leftrightarrow$
$${dn_N-ld_N\over(d_N-n_Ng)(d-lg)}>0
\Leftrightarrow dn_N-ld_N>0$$
ce qui contredit l'hypoth\`ese $\mu(N)={d_N\over
n_N}\geq {d\over l}$. Donc $F$ est
stable.\medskip
\noindent Pour $E$ semi-stable, on
obtient $dn_N-ld_N\geq 0$, mais l'hypoth\`ese de
non-semi-stabilit\'e est ${d_N\over n_N}>{d\over l}$,
ce qui est encore contradictoire.\par \noindent Le
th\'eor\`eme B-1 est enti\`erement d\'emontr\'e.\par
\noindent On a
consid\'er\'e que les espaces
$W_{n,d}^{n+l}$ \'etaient r\'eduit. En fait, c'est le cas, car
on v\'erifie facilement (cf partie C) que le morphisme de
Petri en tout point de  $W_{n,d}^{n+l}$ est injectif.
\hbox to 0,5cm{}$\diamondsuit$\vskip 0,7cm

\noindent{\rmgXIII C- Structure de $\bf W_{n,d}^{k-1}$
pour
$\bf{d\over n}\leq2$}\vskip 0,6cm

\noindent Nous pr\'ecisons la structure des espaces de
Brill-Noether pour les fibr\'es stables de pente $\leq
2$:\bigskip
\noindent{\bf C-1 Proposition:} {\it Soit $d$ et $n$ deux
entiers tels que ${d\over n}\leq 2$. Pour tout entier $k$, si
$W_{n,d}^{k-1}$ est non-vide, alors toutes les
composantes irr\'eductibles de $W_{n,d}^{k-1}$
sont de dimension $\rho(g,n,d,k-1)$ ainsi que
$${\rm
Sing}\,W_{n,d}^{k-1}=W_{n,d}^{k}$$
sauf dans le cas suivant: si la courbe est hyperelliptique,
$d=2n=2$ et $k=2$. Dans ce cas $W_{1,2}^{1}=\{L\}$,
l'unique $g_2^1$ de la courbe.}\medskip
\noindent{\it
D\'emonstration:} Soit $E$ un fibr\'e stable de rang $n$,
de degr\'e $d$ et de pente $\mu(E)\leq2$. On suppose
que  si la courbe est hyperelliptique alors  $E\not
\simeq L$, l'unique $g_2^1$ de la courbe. Il faut montrer
que le morphisme de Petri $$p:{\rm H}^0(E)\otimes {\rm
H}^0(E^*\otimes K) \rightarrow {\rm H}^0(E\otimes
E^*\otimes K)$$ est injectif.\par \noindent On a une suite
exacte $$0\hfl{}{} D(E)^*\hfl{}{} {\rm H}^0(E)\otimes {\cal
O} \hfl{}{} E\hfl{}{} N\hfl{}{}0$$
o\`u $N$ est le faisceau conoyau. On tensorise cette
suite exacte par $E^*\otimes K$
$$0\hfl{}{} D(E)^*\otimes E^*\otimes K
\hfl{}{} {\rm H}^0(E)\otimes E^*\otimes K$$
$$\hfl{}{} E\otimes E^*\otimes K\hfl{}{}
N\otimes E^*\otimes K\hfl{}{}0$$
et en globalisant on obtient que le noyau du morphisme
de Petri est isomorphe \`a
${\rm H}^0(D(E)^*\otimes E^*\otimes K)$, ou
encore d'apr\`es le th\'eor\`eme de dualit\'e de Serre, \`a 
${\rm H}^1(E\otimes D(E))^*$. On veut donc montrer que
${\rm h}^1(E\otimes D(E))=0$. Si $E\simeq E_K$, c'est
clair. Supposons maintenant que $E\not\simeq
E_K$.\medskip
\noindent Pour terminer la d\'emonstration de la
proposition, on montre que le gradu\'e de $D(E)$ n'est
constitu\'e que de fibr\'es semi-stables de pente $\geq
2g$; c'est ce que fait le lemme ci-dessous.\hbox
to 1cm{} $\diamondsuit$\medskip
\noindent{\bf C-2
Lemme:} {\it Soit $E$ un fibr\'e stable de pente
$\mu(E)\leq2$. Le gradu\'e de $D(E)$ contient un fibr\'e
de pente $<2g$ si et seulement si l'on est dans l'un des
cas suivant:\par - $E\simeq E_K$
(la courbe est non-hyperelliptique);\par
- $E\simeq L$ et
la courbe est hyperelliptique de $g_2^1$
$L$.}\medskip
\noindent{\it D\'emonstration du
lemme:}  on suppose  que
$E\not\simeq E_K$ ($C$ non-hyperelliptique) et
$E\not\simeq L$ ($C$ hyperelliptique). Soit $Im$, le
faisceau image du morphisme d'\'evaluation de $E$. On
peut \'ecrire $Im\simeq Im'\oplus {\cal O}^a$, avec
h$^0(Im'^*)=0$. On a alors $D(E)\simeq D(Im')$ et donc 
une suite exacte $$0\hfl{}{} D(E)^*\hfl{}{}{\rm
H}^0(Im')\otimes {\cal O} \hfl{}{} Im'\hfl{}{}0$$ que l'on
dualise  $$0\hfl{}{} Im'^*\hfl{}{}{\rm H}^0(Im')^*\otimes
{\cal O} \hfl{}{} D(E)\hfl{}{}0\ .$$ Cette derni\`ere suite
exacte s'imbrique dans un diagramme commutatif
$$\diagram{ &&0&&0&&0&&&\cr
&&\downarrow&&\downarrow&&\downarrow&&&\cr
0&\longrightarrow&Im'^*&\longrightarrow&  {\rm
H}^0(Im')^*\otimes {\cal O}&\longrightarrow &
D(E)&\longrightarrow&0\cr
&&\downarrow&&\downarrow&&\downarrow&&&\cr
0&\longrightarrow& D\bigl(D(E)\bigr)^*&\longrightarrow&
{\rm H}^0(D(E))\otimes {\cal O}& \longrightarrow&
D(E)&\longrightarrow&0\cr
&&\downarrow&&\downarrow&&\downarrow&&&\cr
0&\longrightarrow&{\cal O}^b&\longrightarrow&{\cal
O}^b  &\longrightarrow &0&&\cr
&&\downarrow&&\downarrow&&&&&\cr
&&0&&0&&&&&\cr }$$
o\`u ${\cal O}^b$ est le conoyau du morphisme
${\rm H}^0(Im')^*\otimes{\cal O}\hookrightarrow 
{\rm H}^0(D(E))\otimes {\cal O}$.
On dualise la suite exacte verticale de gauche:
$$0\hfl{}{} {\cal O}^b\hfl{}{} D\bigl(D(E)\bigr)
\hfl{}{}Im'\hfl{}{}0$$
Notons que $Im'$ est un sous-faisceau de $E$
engendr\'e par ses sections tel que h$^0(Im'^*)=0$.
Donc soit le gradu\'e de $Im'$ n'est compos\'e
que de fibr\'es semi-stables de pente strictement
comprise entre $1$ et $2$ soit $Im'\simeq E$. Dans les
deux cas tous les sous-faisceaux propres de $Im'$ sont
de pente $<2$. La suite exacte ci-dessus montre qu'il
en est de m\^eme pour $D\bigl(D(E)\bigr)$. D'autre part
la proposition A-2 du chapitre 2 implique que
$D(E)\simeq D(Im')$ est de pente $\geq 2g$.\par
\noindent Supposons
maintenant que dans le gradu\'e de $D(E)$, Grad$\,
D(E)=\oplus_{1\leq i\leq p} Q_i$, on ait des fibr\'es $Q_j$
de pente $<2g$. On note $j_0$ le plus petit des $j$
correspondant. Comme le fibr\'e $D(E)$ est de pente
$\geq2g$, $j_0$ est sup\'erieur ou \'egal \`a $2$.   $D(E)$
s'inscrit dans une suite exacte $$0\hfl{}{} H\hfl{}{}
D(E)\hfl{}{} G\hfl{}{}0$$ o\`u Grad$\,H=\oplus_{1\leq
i<j_0}Q_i$ et Grad$\, G=\oplus_{j_0\leq i\leq p}Q_i$. Le
fibr\'e $H$ est engendr\'e par ses sections et 
h$^1(H)=0$, on obtient alors un isomorphisme $${\rm
H}^0(D(E))\simeq {\rm H}^0(H)\oplus {\rm H}^0(G)$$
et donc un diagramme commutatif
$$\diagram{
&&0&&0&&0&&&\cr
&&\downarrow&&\downarrow&&\downarrow&&&\cr
0&\longrightarrow&D(H)^*&\longrightarrow& 
{\rm H}^0(H)\otimes {\cal O}&\longrightarrow
& H&\longrightarrow&0\cr
&&\downarrow&&\downarrow&&\downarrow&&&\cr
0&\longrightarrow& D\bigl(D(E)\bigr)^*&\longrightarrow&
{\rm H}^0(D(E))\otimes {\cal O}&
\longrightarrow& D(E)&\longrightarrow&0\cr
&&\downarrow&&\downarrow&&\downarrow&&&\cr
0&\longrightarrow&D(G)^*&\longrightarrow&
{\rm H}^0(G)\otimes {\cal O}  &\longrightarrow
&G&\longrightarrow&0\cr
&&\downarrow&&\downarrow&&\downarrow&&&\cr
&&0&&0&&0&&&\cr
}$$
dont on d\'eduit que $D(G)$ est un sous-fibr\'e propre de
$D\bigl(D(E)\bigr)$. Le th\'eor\`eme de Clifford
implique que $${\rm h}^0(G)\leq n_G+{d_G\over 2}$$
et donc que
$$\mu(D(G))\geq d_G/{d_G\over 2}=2\ ,$$
mais  tous les sous-faisceaux propres de
$D\bigl(D(E)\bigr)$ sont de pente $<2$. Ce qui
termine la d\'emonstration du lemme.\hbox to
1cm{}$\diamondsuit$\bigskip \noindent {\bf C-3
Remarque:} On devrait en fait pouvoir prouver que si ces
espaces de Brill-Noether sont non-vides, alors ils sont
irr\'eductibles.\vskip 1cm

\noindent{\rmgXV Chapitre 3: Fibr\'es vectoriels
de pente $\bf <2$  }\vskip 0,8cm 
\noindent Nous avons vu que au-dessus de la
droite $\Delta$, il n'existe pas de fibr\'es stables de
pente $<2$ (Chap. 2 th\'eor\`eme A-1). De plus,
pour une pente $1<\mu<2$,  un fibr\'e stable
correspondant
\`a un point sur la droite $\Delta$ est de rang $n$ et
de degr\'e
$d=n+gl$ et poss\`ede $n+l$ sections, pour un
entier
$l$. Le th\'eor\`eme B-1 du Chap.2 montre qu'il
existe de tels fibr\'es. Il est maintenant naturel de se
poser le probl\`eme de l'existence pour des fibr\'es
stables de rang
$n$ et de degr\'e
$d=n+gl +l'$ avec $0<l'<g$ poss\'edant $k$
sections globales ind\'ependantes. Ces fibr\'es ont
en fait au plus
$n+l$ sections. Nous proc\'ederons en deux
\'etapes: on traite d'abord le cas o\`u $l=0$,
c'est-\`a-dire $d<n+g$ (partie A), puis le cas
o\`u $l>0$ (partie B). On aura ainsi trait\'e
tous les cas possibles pour $1<\mu<2$.\vskip 0,7cm

\noindent{\rmgXIII A- Fibr\'es vectoriels avec peu
de sections. }\vskip 0,6cm 

\noindent Dans cette partie, on consid\`ere que
$d<n+g$ et $1<{d\over n}<2$. On montre alors
 qu'il existe des fibr\'es stables de rang
$n$, de degr\'e $d$ poss\'edant $n$ sections
globales et on donne une description des 
espaces de Brill-Noether associ\'es:\bigskip

\noindent{\bf A-1 Th\'eror\`eme:} {\it Si $n$, $d$ et
$k$ v\'erifient
$n<d<n+g$,  ${d \over n}<2$ et $k\leq n$, alors
${
W}_{n,d}^{k-1}$ 
est irr\'eductible de la bonne
dimension et le lieu des singularit\'es est celui
attendu. De plus, ${
W}_{n,d}^{k-1}$ est dense dans ${\widetilde
W}_{n,d}^{k-1}$ . }\medskip

\noindent{\it D\'emonstration:} La d\'emonstration
reprend exactement les m\^emes arguments que
ceux donn\'es dans [{\bf B-G-N}]. On va montrer
que si $E$ est un fibr\'e semi-stable de rang $n$ et
de degr\'e $d$,
$n$ et $d$ comme ci-dessus et si
h$^0(E)=k$, alors on a une suite exacte de
faisceaux:
$$0\longrightarrow {\cal O}^k\longrightarrow E
\longrightarrow S\longrightarrow 0$$
o\`u $S$ peut avoir de la torsion. Ensuite, le
probl\`eme se ram\`ene \`a compter ces
extensions. Le lemme ci-dessous montre
l'existence de telles suites exactes:\bigskip

\noindent{\bf A-2 Lemme:} {\it Soit $E$ un fibr\'e
vectoriel semi-stable de rang $n$, de degr\'e $d$
tel que $d<n+g$ et ${d\over n}<2$. Soit $V\subset
{\rm H}^0(E)$, un espace vectoriel de sections
globales. Alors on a un morphisme injectif de
faisceaux:\vskip -0.2cm
$$V\otimes {\cal O}\hookrightarrow E$$ }\medskip

\noindent{\it D\'emonstration du lemme}: C'est une
cons\'equence direct du th\'eor\`eme A-1 du
chapitre 2. En effet, si $M$,
 le faisceau image du morphisme
$V\otimes {\cal O}\rightarrow E$, n'est pas le
faisceau trivial, alors il admet un sous-faisceau 
$M'$ tel  que h$^0(M'^*)=0$ et comme $M'$ est un
sous-faisceau de $E$, on a
$\mu(M')={d'\over n'}<2$. Le th\'eor\`eme A-1 du
chap. 1 dit alors que
 $${\rm h}^0(M')\leq n'+{1\over g}(d'-n')$$
\noindent Par hypoth\`ese h$^0(M')>n'$, donc
${1\over g}(d'-n')\geq 1$, mais alors ${d'\over
n'}\geq1+{g\over n'}$ et ceci contredit la
semi-stabilit\'e de $E$: 
$\mu(M')\geq 1+{g\over n'}\geq 1+{g\over
n}>\mu(E)$. Ce qui termine la d\'emonstration du
lemme.\hbox to0,5cm{}$\diamondsuit$\medskip

\noindent Notons que si $S$ est un faisceau
coh\'erent, comme on est sur une courbe, $S$ peut
s'\'ecrire
$S\simeq G\oplus T$ o\`u $G$ est un fibr\'e
vectoriel et $T$ est un faisceau de torsion (c'est le
th\'eor\`eme de d\'ecomposition d'un module de
type fini au-dessus d'un id\'eal principal).

\noindent On va donc \'etudier les extensions
$$0\longrightarrow {\cal O}^k\longrightarrow E
\longrightarrow G\oplus T\longrightarrow 0$$
o\`u:\par
- $E$ est un fibr\'e stable de rang $n$ et de degr\'e
$d$.\par
- $G$ est un fibr\'e vectoriel avec 
h$^0(G^*)=0 $
car $G$ est un quotient de $E$.\par
- $T$ est un faisceau de points de degr\'e
$\delta$.\medskip

\noindent Ces extensions sont
param\'etr\'ees par Ext$^1(G\oplus T,{\cal O}^k)
\simeq
k.\hbox{Ext}^1(G\oplus T,{\cal O})$.\par

\noindent Pour calculer la dimension de
$\hbox{Ext}^1(G\oplus T,{\cal O})$, on 
utilise la dualit\'e de Serre

$$\hbox{Ext}^1(G\oplus T,{\cal O})\simeq
{\rm H}^0\bigl((G\oplus T)\otimes K)^*
\simeq {\rm H}^1(G^*)\oplus  {\rm H}^0(T)^*$$

\noindent Comme par hypoth\`ese h$^0(G^*)=0$,
et que h$^0(T)=\delta$,
on en d\'eduit que
$$\hbox{dim}\,\hbox{Ext}^1(G\oplus T,{\cal
O})=d+(n-k)(g-1)\ .$$
Donc ces extensions sont 
param\'etr\'ees par des $k$-uplets $(e_1,\cdots,e_k)$,
$e_i\in
\hbox{Ext}^1(G\oplus T,{\cal O})$ et deux de ces
extensions donnent deux fibr\'es isomorphes si les
$k$-uplets correspondants sont dans la 
m\^eme orbite par
l'action naturelle du groupe GL($k$). 
On en d\'eduit que si
$(e_1,\cdots,e_k)$ sont lin\'eairement d\'ependants,
on peut supposer que $e_k=0$, mais alors ${\cal O}$
est un facteur direct de $E$, ce qui contredit la
stabilit\'e de $E$. On remarque que 
$(e_1,\cdots,e_k)$
sont n\'ecessairement d\'ependants si
$k>d+(n-k)(g-1)$ ou encore si $n>d+(n-k)g$. Ce
qui ne peut arriver que si $n>d$, c.-\`a-d. si $\mu<1$,
cas trait\'e dans  {\bf [B-G-N]}. \par
\noindent Les extensions telles que
$E$ n'admettent pas ${\cal O}$ comme quotient,
sont
donc classifi\'ees par la grassmanienne des
sous-espaces vectoriels de dimension $k$ de
Ext$^1(G\oplus T,{\cal O})$. On a:
$$\hbox{dim}\biggl(\hbox{Grass}_k\bigl(\hbox{Ext}^1
(G\oplus T,{\cal O})\bigr)\biggr)=k(d+(n-k)g-n)\ .$$
Nous utiliserons cette grassmanienne pour la
question de l'existence, mais pour le calcul de
dimension, il faut diminuer le nombre de
param\`etres dont d\'ependent ces extensions. 
Il est
alors n\'ecessaire de pr\'eciser l'action du groupe
\hbox{Aut$({\cal O}^k)\times\hbox{Aut}(G\oplus T)$}
sur Ext$^1(G\oplus T,{\cal O}^k)$:\bigskip

\noindent Soit
$$0\longrightarrow {\cal O}^k\longrightarrow E
\longrightarrow G\oplus T\longrightarrow 0$$
une  extension avec $E$ stable. Alors
le stabilisateur de cette action en ce point est
l'ensemble $\{(\lambda \hbox{Id}_{{\cal
O}^k},1/\lambda
\hbox{Id}_{G\oplus T})\mid\ \lambda\in\hbox{\mathsy
C}\}$, donc il est de dimension $1$.\par

 \noindent En effet, soit $(l,h)\in
\hbox{Aut}({\cal O}^k)\times\hbox{Aut}(G\oplus T)$.
L'action de $(l,h)$ sur
l'extension\par
$$0\longrightarrow {\cal O}^k\hfl{i}{} E
\hfl{j}{} G\oplus T\longrightarrow 0$$
donne l'extension
$$0\longrightarrow {\cal O}^k\hfl{i\circ l}{} E
\hfl{j\circ h}{} G\oplus T
\longrightarrow0\  .$$
Ces deux extensions repr\'esentent le m\^eme
\'el\'ement dans Ext$^1(G\oplus T,{\cal O}^k)$ 
si et seulement si il existe
un diagramme commutatif:
$$\matrix{
0&\longrightarrow &{\cal O}^k&\hfl{i}{} &E&\hfl{j}{} 
&G\oplus T&
\longrightarrow&0\cr
&&\parallel&&\vfl{\displaystyle\wr\kern-2pt}{\varphi}
&&\parallel&&\cr
0&\longrightarrow& {\cal O}^k&\hfl{i\circ l}{} &
E&\hfl{j\circ h}{} 
&G\oplus T&\longrightarrow&0\cr
}$$
$\varphi$ \'etant un automorphisme de $E$, donc
$\varphi=\lambda\hbox{Id}_E$ par hypoth\`ese.
Au niveau
des fibres en un point $x$, on obtient le diagramme
commutatif de ${\cal O}_x$-modules:
$$\matrix{
0&\longrightarrow &{\cal O}^k_x&\hfl{i_x}{} &E_x&
\hfl{j_x}{} &
G_x\oplus T_x&\longrightarrow&0\cr
&&\parallel&&\vfl{\displaystyle\wr
\kern-2pt}{\lambda.Id_x}&&\parallel&&\cr 
0&\longrightarrow&
{\cal O}^k_x&\hfl{i_x\circ l_x}{} & E_x&
\hfl{j_x\circ h_x}{}
&G_x\oplus T_x&\longrightarrow&0\cr
}$$
et il est clair que $l_x=\lambda\hbox{Id}_{{\cal
O}^k_x}$ et que
$h_x=1/\lambda \hbox{Id}_{{(G\oplus T)}_x}$, d'o\`u
le r\'esultat.
\medskip
\noindent Donc si $U$ est l'ouvert de
Ext$^1(G\oplus T,{\cal O}^k)$ correspondant aux
$E$ stables, la propri\'et\'e universelle de la
vari\'et\'e de modules $U_{n,d}$ implique
l'existence d'un morphisme
$$U\hfl{}{} U_{n,d}$$
dont les fibres sont de dimension sup\'erieure
\`a
$$\hbox{dim Aut}\,{\cal O}^k+
\hbox{dim Aut}\,(G\oplus T)-1\ .$$
\noindent Pour $G$ et $T$ fix\'es, on en d\'eduit
 que l'image de ce morphisme d\'epend de au plus
 
$$\chi=k(d+(n-k)g-n)-\hbox{dim}\,\hbox{Aut}(G\oplus
T)+1$$ param\`etres, si elle est non vide.
On va minorer 
$\hbox{dim}\,\hbox{Aut}(G\oplus T)=m$:\medskip
- Si $T=0$, dim Hom$(G,G)\geq 1$, 
donc $m\geq 1$ et
$\chi\leq k(d+(n-k)g-n)$\par
- Si $T,G\not=0$, on a: Hom$(T,G)$=0, dim
Hom$(G,G)\geq1$, dim Hom$(G,T)=\delta(n-k)$
et les
automorphismes de $T$ sont r\'eduits 
aux homoth\'ethies
en chaque point, c.-\`a-d. dim Aut$\,T\geq$
nombre de
points du support de $T$.
Donc $m\geq \delta (n-k)+2$ et
$\chi\leq k(d+(n-k)g-n)-\delta (n-k)-1$\par
- Si $G =0$, alors $k=n$ et $m\geq$ nombre de
points du support de $T$, et donc $\chi\leq
n(d-n)-m+1$.\medskip

\noindent On peut alors montrer que l'espace de
Brill-Noether,
$W_{n,d}^{k-1}$, est irr\'eductible, de la bonne
dimension avec un bon lieu de singularit\'e s'il
est non vide.\par
\noindent En effet, le lemme A-2 montre que tout
fibr\'e stable $E$
poss\'edant $k$ sections globales
ind\'ependantes s'\'ecrit comme une extension
$$0\longrightarrow {\cal O}^k\longrightarrow 
E\longrightarrow
G\oplus T\longrightarrow0\ . $$

\noindent Supposons $k<n$, c'est-\`a-dire
$G\not=0$. Pour $G$ et  $T$ fix\'es, ces
extensions d\'ependent de au plus 
$k(d+(n-k)g-n)$ param\`etres,
$\bigl(\,k(d+(n-k)g-n)-\delta (n-k)-1$ param\`etres si
$T\not =0\bigr)$,
 et comme 
${\rm H}^0(G^*)=0$ pour $d$, $k$ et $T$ fix\'es
l'ensemble des $G$
consid\'er\'es forment une famille limit\'ee,
cette famille
d\'epend de au plus $(n-k)^2(g-1)+1$ 
param\`etres (cf {\bf [B-G-N]} lemme 4.1). Le
choix de $T$ pour $\delta$ fix\'e rajoute $\delta$
param\`etres donc on obtient la majoration
$$\matrix{
\hbox{dim}\,W_{n,d}^{k-1}&\leq &(n-k)^2(g-1)+1+
k(d+(n-k)g-n)\cr
&\leq&\rho(g,d,n,k-1)\hbox to 4,5cm{}\cr}$$
et on en d\'eduit que si $W_{n,d}^{k-1}$ est non vide
alors il est de la bonne dimension.\par
\noindent De plus on remarque que si $T\not=0$, 
alors les
extensions d\'ependent de au plus
$\rho(g,d,n,k-1)-1$
param\`etres. Or on sait que toute composante
irr\'eductible de $W_{n,d}^{k-1}$ est de dimension au
moins $\rho(g,d,n,k-1)$ donc toute composante de
$W_{n,d}^{k-1}$ contient un sous-ensemble dense
${\cal V}$ contenant  les extensions
\centerline {$0\longrightarrow {\cal
O}^k\longrightarrow E\longrightarrow
G\longrightarrow0\ .$}\medskip 
\noindent Ces extensions
sont param\'etr\'ees par une vari\'et\'e projective
irr\'eductible (cf  {\bf [N-R]} proposition 2.6 ou
 [{\bf B-G-N}]). L'ouvert de cette vari\'et\'e
correspondant aux fibr\'es stables s'envoie
surjectivement sur ${\cal V}\subset W_{n,d}^{k-1}$.
On en d\'eduit l'irr\'eductibilit\'e de
$W_{n,d}^{k-1}$. La proposition  C-1 du Chap. 2 montre
montre l'assertion sur les points singuliers  de
$W_{n,d}^{k-1}$.\medskip

\noindent Si $k=n$, alors $G=0$, 
$\delta=d$ et on peut
fixer le nombre $p$ de points du
support de $T$ et le
degr\'e en chacun de ces points. Ces valeurs sont
discr\`etes. Alors, d'apr\`es les calculs pr\'ec\'edents,
les suites exactes
$$0\longrightarrow {\cal O}^k \longrightarrow E
\longrightarrow T\longrightarrow 0$$
d\'ependent de au plus $n(d-n)-p+1$
param\`etres et le
choix de $T$ ajoute $p$ param\`etres.
On retrouve bien
$$\hbox{dim}\,W_{n,d}^{n-1}\leq n(d-n)+1=\rho\ .$$
 Pour l'irr\'eductibilit\'e,
les arguments sont encore les m\^emes:
Quot$^d(n,C)$, la vari\'et\'e param\'etrisant les
faisceaux de torsion de degr\'e $d$ quotient de ${\cal
O}^n$, est irr\'eductible (cf [{\bf Rego}]); on en d\'eduit
l'existence d'une vari\'et\'e irr\'eductible qui param\'etrise
toutes les extensions ci-dessus et l'ouvert correspondant aux
fibr\'es stables s'envoie surjectivement sur
$W_{n,d}^{n-1}$. Et on termine avec la proposition
C du Chap. 2.\medskip

\noindent Pour l'existence, on adapte encore
des id\'ees qui se trouvent dans {\bf [BGN]}.\par
\noindent On \'etudie les extensions
$$0\longrightarrow {\cal O}^n\longrightarrow E
\longrightarrow \Theta\longrightarrow0$$
o\`u: \par
- $\Theta$ est un faisceau de torsion
de degr\'e $d$: par exemple $\Theta$ est de
support $d$ points distincts.\par 
- $E$ est un fibr\'e
tel que ${\cal O}$ ne soit pas un facteur
direct de $E$.\medskip

\noindent D'apr\`es ce qui pr\'ec\`ede, il
existe  de tels fibr\'es et
ces extensions sont param\'etr\'ees
par la grassmanienne
$$\hbox{Grass}_n\hbox{Ext}^1(\Theta,{\cal O})$$
qui est une vari\'et\'e de dimension $n(d-n)>0$.\par
\noindent Si $E$ est non stable, alors il existe un 
fibr\'e
stable $F$ quotient de $E$ tel que 
$\mu(F)\leq \mu(E)$. On obtient un diagramme
commutatif
$$\matrix{
&&0&&0&&0&&&\cr
&&\downarrow&&\downarrow&&\downarrow&&&\cr
0&\longrightarrow& M'&\longrightarrow& F'&
\longrightarrow& \Theta''&\longrightarrow&0\cr
&&\downarrow&&\downarrow&&\downarrow&&&\cr
0&\longrightarrow& {\cal O}^n&\longrightarrow& E&
\longrightarrow& \Theta&\longrightarrow&0\cr
&&\downarrow&&\downarrow&&\downarrow&&&\cr
0&\longrightarrow& M&\longrightarrow& F&
\longrightarrow& \Theta'&\longrightarrow&0\cr
&&\downarrow&&\downarrow&&\downarrow&&&\cr
&&0&&0&&0&&&\cr
}$$
o\`u $M$ est l'image de ${\cal O}^n\rightarrow F$,
$M'$ et $F'$ sont les noyaux des morphismes
verticaux et $\Theta'$ et $\Theta''$ sont
des faisceaux de torsions de degr\'e $d'$ et
$d-d'$.\par 
\noindent $F$ \'etant stable et 
$\mu(F)\leq\mu(E)=1+{g\over n}$, le lemme A-2
implique que l'on a 
encore un morphisme injectif
$${\rm H}^0(F)\otimes{\cal O}\hookrightarrow F\ .$$
$M$ est un sous-faisceau de $F$
engendr\'e par ses sections, on doit donc avoir
$M\simeq {\cal O}^l$. Notons que $F\not \simeq {\cal
O}$ puisque ${\cal O}$ n'est pas un facteur direct de
$E$. \par

\noindent Le diagramme
devient: 
$$\matrix{ &&0&&0&&0&&&\cr
&&\downarrow&&\downarrow&&\downarrow&&&\cr
0&\longrightarrow&{\cal O}^{n-l}&\longrightarrow& F'&
\longrightarrow& \Theta''&\longrightarrow&0\cr
&&\downarrow&&\downarrow&&\downarrow&&&\cr
0&\longrightarrow& {\cal O}^n&\longrightarrow& E&
\longrightarrow& \Theta&\longrightarrow&0\cr
&&\downarrow&&\downarrow&&\downarrow&&&\cr
0&\longrightarrow& {\cal O}^l&\longrightarrow& F&
\longrightarrow& \Theta'&\longrightarrow&0\cr
&&\downarrow&&\downarrow&&\downarrow&&&\cr
&&0&&0&&0&&&\cr
}$$
Comme dans {\bf [BGN]}, on a un morphisme
surjectif $$\hbox{Ext}^1(\Theta,{\cal O})\rightarrow 
\hbox{Ext}^1(\Theta'',{\cal O})$$ et l'existence
de la suite exacte horizontale sup\'erieure implique
que l'image du n-uplet $(e_1,\cdots,e_n)\in
\hbox{Ext}^1(\Theta,{\cal O})$ repr\'esentant
l'extension
$$0\hfl{}{}{\cal O}^n\hfl{}{} E\hfl{}{} \Theta\hfl{}{}0$$
a pour image par ce morphisme un n-uplet dont
au plus $n-l$ \'el\'ements sont ind\'ependants. 
Ceci d\'etermine une sous-vari\'et\'e $Z$ de
$\oplus^n \hbox{Ext}^1(\Theta,{\cal O})$ de
codimension
$l(d-d'-n+l)$.\par
\noindent Or, par hypoth\`ese, $d>n$ et 
${d\over n}\geq {d'\over
l}$. De plus, la suite exacte
horizontale du bas montre que $F$ est un fibr\'e
stable de rang $l$ poss\'edant $l$ sections globales
ind\'ependantes donc que $d'>l$ (cf th\'eor\`eme A-1
du Chap. 2).  Donc,  $$d-d'\geq ({n\over
l}-1)d'=(n-l){d'\over l}>(n-l)\ .$$ On en d\'eduit que
$$\hbox{codim}\,Z\geq l$$ Il existe donc des extensions
$$0\hfl{}{}{\cal O}^n\hfl{}{} E\hfl{}{} \Theta\hfl{}{}0$$
qui ne peuvent se compl\'eter en un diagramme
comme ci-dessus pour aucune valeur de $d'$ et le
fibr\'e
$E$ correspondant est stable et poss\`ede $n$
sections globales ind\'ependantes. Le th\'eor\`eme
A-1 est ainsi enti\`erement d\'emontr\'e.\hbox to
0,5cm{}
$\diamondsuit$\vskip 0,7cm

\noindent{\rmgXIII B- Fibr\'es vectoriels avec
beaucoup de sections.}\vskip 0,6cm
\noindent Ici, on consid\`ere que $d=n+gl+l'$, avec
${d\over n}<2$, $0<l'<g$ et $l>0$. L'id\'ee est de d\'eduire
l'existence de fibr\'es $E$, de rang $n$ et de degr\'e
$d$, de l'existence de fibr\'es stables
correspondant \`a des points sur la droite $\Delta$. Notons
tout d'abord que le th\'eor\`eme A-1 du Chap. 2 implique que
$${\rm h}^0(E)\leq n+{1\over g}(d-n)$$
et comme $d=n+gl +l'$, on en d\'eduit que
h$^0(E)\leq n+l$. On va donc chercher des fibr\'es
stables $E$ tel que h$^0(E)=n+l$. On
d\'eduira de l'existence de ces fibr\'es $E$
le th\'eor\`eme ci-dessous, qui termine
l'\'etude des fibr\'es vectoriels stables de
pente $<2$:\bigskip

\noindent{\bf B-1 Th\'eor\`eme:} {\it On suppose
que
$d=n+gl+l'
$ o\`u $n$, $l$, $l'$ sont des entiers v\'erifiant 
${d\over n}<2$, $ 0<l'<g$ et $l>0$. Alors
$W_{n,d}^{n+l-1}$ est non vide.}\medskip

\noindent{\it D\'emonstration:} Le th\'eor\`eme B-1
du Chap. 2 donne l'existence de fibr\'es stables $E'$
de rang
$m=n+l'$ et de degr\'e $d=n+gl+l'$
 tel que $E'$ s'inscrit dans une suite exacte

$$0\hfl{}{}D(E')^*\hfl{}{} {\cal O}^{m+l}\hfl{}{}
E'\hfl{}{}0$$
avec $D(E')$ stable.\par
\noindent On va en fait chercher les fibr\'es $E$
comme des quotients de $E'$ ($ l'<n$):
$$0\hfl{}{}{\cal O}^{l'}\hfl{}{}E' \hfl{}{}
E\hfl{}{}0$$
Il est clair que si $E$ est stable, alors $E$ v\'erifie
toutes les propri\'et\'es voulues.\medskip

\noindent Nous utiliserons la terminologie suivante: nous
avons montr\'e que l'on a un isomorphisme entre
$W_{m,m+gl}^{m+l-1}$ et
$U_{l,m+gl}$. Cet isomorphisme \'etant d\'efini,
avec les notations ci-dessus, par
$E'\mapsto D(E')$. On dira que $E'$ est g\'en\'erique
dans $W_{m,m+gl}^{m+l-1}$ si $D(E')$ est
g\'en\'erique dans le sens habituel dans
$U_{l,m+gl}$. \medskip

\noindent 
Soit une extension
$$0\hfl{}{}{\cal O}^{l'}\hfl{}{}E' \hfl{}{}
E\hfl{}{}0\ ,$$
et supposons que $E$ n'est pas stable. Soit $G$
un fibr\'e quotient stable de pente minimale:
$$0\hfl{}{}H\hfl{}{}E \hfl{}{}
G\hfl{}{}0\ .$$
Le fibr\'e $G$ est stable de pente
$\mu(G)={d_G\over n_G}\leq \mu(E)={n+gl+l'\over
n}=1+{gl+l'\over n}<2$. Le th\'eor\`eme A-1 du
Chap. 2 implique que $${\rm h}^0(G)\leq
n_G+{1\over g}(d_G-n_G)\ .$$
En posant
h$^0(G)=n_G+l_G$, on obtient
$$\mu(G)\geq
1+{gl_G\over n_G}\ .$$
Et comme $\mu(G)\leq
\mu(E)=1+{gl+l'\over n}$, on en d\'eduit que
$l_G\leq l$ (on a $l'<g $).\medskip 
\noindent Expliciter directement sur $E$ le fait que
$E$ ne soit pas stable semble fort complexe et ne
pas aboutir. L'id\'ee est d'\'etudier les
cons\'equences de $E$ non stable sur les fibr\'es
$E'$ dont nous avons une bonne description. C'est
ce que nous permet de faire le diagramme
commutatif ci-dessous obtenu \`a partir des deux
suites exactes pr\'ec\'edentes:
$$\matrix{
&&0&&0&&0&&&\cr
&&\downarrow&&\downarrow&&\downarrow&&&\cr
0&\longrightarrow&0&\longrightarrow& 
0&\longrightarrow
& H&&\cr
&&\downarrow&&\downarrow&&\downarrow&&&\cr
0&\longrightarrow& {\cal O}^{l'}&\longrightarrow&
E'&
\longrightarrow& E&\longrightarrow&0\cr
&&\downarrow&&\downarrow&&\downarrow&&&\cr
0&\longrightarrow&Q&\longrightarrow&E' 
&\longrightarrow
&G&\longrightarrow&0\cr
&&\downarrow&&\downarrow&&\downarrow&&&\cr
&&H&\longrightarrow&0&&0&&&\cr
&&\downarrow&&&&&&&\cr
&&0&&&&&&&\cr
}$$
\centerline{\it Diagramme 1}\medskip
\noindent o\`u $Q$ est le noyau du morphisme
$E'\rightarrow G$ et l'apparition du fibr\'e $H$ en
bas \`a gauche est donn\'ee par le lemme
d'alg\`ebre appel\'e "{\it Le diagramme du serpent"}
(cf [{\bf Pe}]).
$Q$ est de pente
$\mu(Q)={d_Q\over n_Q}$.\par
\noindent Pour  d\'emontrer l'existence de fibr\'es
stables $E$ comme annonc\'es, nous allons
d'abord
\'etudier les fibr\'es $Q$, ensuite les compter dans
le cas $l'=1 $. On obtiendra que l'on n'a pas assez
de fibr\'es
$Q$ pour "englober" toutes les sections de
$E'$. Puis on fait une r\'ecurrence sur $l'$.\medskip

\noindent{\bf \'Etude des fibr\'es $\bf Q$:} 
L'objet essentiel de cette \'etude est de montrer
que $Q$ est
g\'en\'eriquement engendr\'e par ses sections. En
fait on ne pourra le montrer directement. On donne
ici tous les \'el\'ements n\'ecessaires et il nous faudra
reprendre cette \'etude quand on aura pos\'e les hypoth\`eses
de la r\'ecurrence. On note h$^0(Q)=n_Q+l_Q$. Comme
$l_G\leq l$, on en d\'eduit que $l_Q\geq 0$. $Q$ a
donc au moins $n_Q$ sections. \par \noindent Tout
d'abord montrons que h$^0(Q^*)=l_Q'< l'$: en effet,
le fibr\'e
$H$ a, par hypoth\`ese, une filtration de
Harder-Narasimhan (cf partie A du Chap. 1) dont
les quotients sont tous de pente
$\geq \mu(G)>0$. Donc h$^0(H^*)=0$. Dans le
diagramme 1, le dual de la suite exacte
verticale de gauche 
$$0\hfl{}{} H^*\hfl{}{} Q^*\hfl{}{}
{\cal O}^{l'}\hfl{}{}0 $$
montre que h$^0(Q^*)\leq l'$. L'existence de $l_Q'$
sections se traduit par un morphisme ${\cal
O}^{l_Q'}\rightarrow Q^* $.  Le sch\'ema ci-dessous
$$\matrix{
&&&& 
{\cal O}^{l_Q'}&
& &&\cr
&&&&\downarrow&&&&&\cr
0&\longrightarrow&H^* &\longrightarrow&
Q^*&
\longrightarrow& {\cal
O}^{l'}&\longrightarrow&0\cr}$$
montre alors que ${\cal O}^{l_Q'}\hookrightarrow
{\cal O}^{l'}$, puisque h$^0(H^*)=0$ et donc que
$Q$ admet ${\cal O}^{l_Q'}$ en facteur direct. On
note $Q\simeq {\cal O}^{l_Q'}\oplus Q'$. Si
h$^0(Q^*)=l'$ alors on obtient 
$Q\simeq {\cal O}^{l'}\oplus H$. Or $Q$ est un
sous-fibr\'e de $E'$ stable et par hypoth\`ese
$\mu(H)\geq \mu(E)>\mu(E')$. $H$ ne peut donc
pas \^etre un sous-fibr\'e de $E'$ et on a bien
h$^0(Q^*)<l'$. \medskip

\noindent  On a deux suites exactes
$$0\hfl{}{} Q'\hfl{}{} Q\hfl{}{}
{\cal O}^{l_Q'}\hfl{}{}0 $$
et
$$0\hfl{}{} Q\hfl{}{} E'\hfl{}{}
G\hfl{}{}0 $$
qui s'imbriquent dans un diagramme commutatif
$$\matrix{
&&0&&0&&0&&&\cr
&&\downarrow&&\downarrow&&\downarrow&&&\cr
0&\longrightarrow&0&\longrightarrow& 
0&\longrightarrow
& {\cal O}^{l_Q'}&&\cr
&&\downarrow&&\downarrow&&\downarrow&&&\cr
0&\longrightarrow& Q'&\longrightarrow&
E'&
\longrightarrow& G'&\longrightarrow&0\cr
&&\downarrow&&\downarrow&&\downarrow&&&\cr
0&\longrightarrow&Q&\longrightarrow&E' 
&\longrightarrow
&G&\longrightarrow&0\cr
&&\downarrow&&\downarrow&&\downarrow&&&\cr
&&{\cal O}^{l_Q'}&\longrightarrow&0&&0&&&\cr
&&\downarrow&&&&&&&\cr
&&0&&&&&&&\cr
}$$
Le fibr\'e $G'$ \'etant le conoyau du morphisme
$Q'\rightarrow E'$. $G'$ est un quotient de $E'$,
donc h$^0(G'^*)=0$. De plus, la suite exacte
verticale de droite montre que le sous-fibr\'e
maximal de $G'$ est de pente $\leq
\mu(G)<2$. D'apr\`es la proposition A-2 du Chap. 1,
on obtient
$${\rm }h^0(G')\leq n_{G'}+{1\over
g}(d_{G'}-n_{G'})\ .$$
On a la m\^eme in\'egalit\'e pour
$Q'$ puisque par construction h$^0(Q'^*)=0$ et que
$Q'$ est un sous-fibr\'e de $E'$:
$${\rm }h^0(Q')\leq n_{Q'}+{1\over
g}(d_{Q'}-n_{Q'})\ .$$
Dans le diagramme ci-dessus, la suite exacte
horizontale du milieu
$$0\hfl{}{} Q'\hfl{}{} E'\hfl{}{}
G'\hfl{}{}0 $$
montre que
$${\rm h}^0(E')\leq {\rm h}^0(G')+{\rm
h}^0(Q')\leq n_{G'}+{1\over
g}(d_{G'}-n_{G'})+n_{Q'}+{1\over
g}(d_{Q'}-n_{Q'})\ .$$
Comme, par hypoth\`ese, h$^0(E')=m+{1\over
g}(d-m)$ (rappelons que rg$\,E'=n+l'=m$),
$m=n_{Q'}+n_{G'}$ et $d=d_{Q'}+ d_{G'}$,
l'in\'egalit\'e ci-dessus implique que
$${\rm }h^0(Q')= n_{Q'}+{1\over
g}(d_{Q'}-n_{Q'})$$
et que
$${\rm }h^0(G')= n_{G'}+{1\over
g}(d_{G'}-n_{G'})\ .$$
Les fibr\'es $G'$ et $Q'$
correspondent \`a des points sur la droite
$\Delta$. Comme h$^0(Q')+l'_Q={\rm h}^0(Q)$ et
$d_{Q'}=d_Q$, on en d\'eduit que h$^0(Q')=n_{Q'}+l_Q$ et
que $d_Q=n_Q-l'_Q+gl_Q$. De la m\^eme fa\c con, on
obtient h$^0(G')= n_{G'}+l_G$ et
$d_G=n_G+gl_G+l'_Q$.\medskip

\noindent Nous allons proc\`eder de la fa\c con suivante:
on va
donner
$\Gamma$ une majoration du nombre de param\`etres 
dont d\'ependent les fibr\'es $Q$. Chacun de ces fibr\'es $Q$
a $n_Q+l_Q$ sections globales
et donc d\'efinit dans
Grass$_{l'}\bigl({\rm H}^0(E')\bigr)$ un
sous-espace de dimension $l'(n_Q+l_Q-l')$, 
la dimension de la
grassmanienne Grass$_{l'}\bigl({\rm H}^0(Q)\bigr)$.\par
\noindent Le sous-espace de Grass$_{l'}\bigl({\rm
H}^0(E')\bigr)$ d\'efini par tous les fibr\'es $Q$ 
d\'epend alors de au plus
$$\Gamma+l'(n_Q+l_Q-l')$$
param\`etres.\par
\noindent Or, la grassmanienne Grass$_{l'}\bigl({\rm
H}^0(E')\bigr)$ est de dimension \hbox{$l'(n_Q+n_G+l-l')$},
donc si l'on montre que
$$l'(n_Q+n_G+l-l')>\Gamma+l'(n_Q+l_Q-l')\ ,$$
alors il existe des sous-espaces ${\cal O}^{l'}\subset {\rm
H}^0(E')$ qui ne s'injectent pas dans un fibr\'e $Q$ comme
ci-dessus et donc le conoyau de l'injection ${\cal O}^{l'}
\hookrightarrow E'$ est g\'en\'eriquement un fibr\'e stable, ce
que nous voudrions. Comme $l=l_Q+l_G$, l'in\'egalit\'e
ci-dessus se simplifie:
$$l'(n_G+l_G)>\Gamma\ \ 
\ \ \ \ \ \ \ (*)$$
Il faut donc trouver une bonne majoration $\Gamma$ du
nombre de param\`etres dont d\'ependent les fibr\'es $Q$.
Pour cela nous aurons besoin du fait que $Q$ est
g\'en\'eriquement engendr\'e par ses sections.

\noindent {\bf Si $\bf l'=1$}, la premi\`ere \'etape de la
r\'ecurrence, alors h$^0(Q^*)<l'$ donne
h$^0(Q^*)=0$ et donc 
 $Q=Q'$ et $G=G'$ et\par
 - $d_Q=n_Q+gl_Q$;
$d_G=n_G+gl_G$.\par
- h$^0(Q)=n_Q+l_Q$;
h$^0(G)=n_G+l_G$.\par
\noindent De plus, dans ce cas, $Q$ s'inscrit dans une
suite exacte (cf diagramme 1)
$$0\hfl{}{}{\cal O}\hfl{}{} Q\hfl{}{} H\hfl{}{}0$$
o\`u $H$ contredit la stabilit\'e de $E$,
c'est-\`a-dire $\mu(H)\geq \mu(E)>\mu(E')$. Ceci
implique que $\mu(H)={d_Q\over
n_Q-1}=1+ {gl_Q+1\over n_Q-1}>\mu(E')$.\par
\noindent Soit
${\cal F}$, le sous-faisceau de $Q$ engendr\'e par
ses sections. On a h$^0({\cal F})=n_Q+l_Q$. Si le
rang de ${\cal F}$ est strictement inf\'erieur \`a
$n_Q$, on d\'eduit de la proposition A-2 du Chap.
2 que ${\cal F}$ a un sous-faisceau ${\cal F}'$ de
pente
$\geq 1+{g(l_Q+1)\over {\rm rg}\,{\cal F}}$. ${\cal
F}'$ serait alors un sous-faisceau de $E'$ stable
de pente $\geq 1+{g(l_Q+1)\over {\rm rg}\,{\cal F}}
>\mu(H)>\mu(E')$, ce qui est absurde. Donc ${\cal
F}$ est de rang $n_Q$. On obtient bien que
$Q$ est g\'en\'eriquement engendr\'e par ses
sections.\medskip
\noindent Pour \'eviter de faire deux fois les m\^emes calculs,
on va supposer que pour
$l'<g$ quelconque et pour 
$E'$ et ${\cal O}^{l'}$ assez g\'en\'eraux, alors $Q$ est
g\'en\'eriquement engendr\'e par ses sections (ceci sera
montr\'e avec la r\'ecurrence). \par

\noindent On suppose donc que le morphisme d'\'evaluation
de
$Q$ s'inscrit dans une suite
exacte
$$0\hfl{}{}D(Q)^*\hfl{}{} {\rm H}^0(Q)\otimes {\cal
O}\hfl{}{} Q\hfl{}{} \Theta\hfl{}{} O $$
o\`u $\Theta $ est le faisceau de torsion
conoyau. On note $d_{\Theta}$ son degr\'e.\par
\noindent On obtient alors un diagramme
commutatif:
$$\matrix{
&&0&&0&&0&&&\cr
&&\downarrow&&\downarrow&&\downarrow&&&\cr
0&\longrightarrow&D(Q)^*&\longrightarrow& 
D(E')^*&\longrightarrow
&D(G)^*&&\cr
&&\downarrow&&\downarrow&&\downarrow&&&\cr
0&\longrightarrow&  {\rm H}^0(Q)\otimes {\cal
O}&
\longrightarrow&
{\rm H}^0(E')\otimes {\cal
O}&
\longrightarrow& {\rm H}^0(G)\otimes {\cal
O}&\longrightarrow&0\cr
&&\downarrow&&\downarrow&&\downarrow&&&\cr
0&\longrightarrow&Q
&\longrightarrow&E' 
&\longrightarrow &G&\longrightarrow&0\cr
&&\downarrow&&\downarrow&&\downarrow&&&\cr
&&\Theta&\longrightarrow&0&&0&&&\cr
&&\downarrow&&&&&&&\cr
&&0&&&&&&&\cr
}$$
\centerline{\it Diagramme 2}\par

\noindent Consid\`erons la suite exacte du "diagramme du
serpent"
$$0\hfl{}{} D(Q)^*\hfl{}{} D(E')^*\hfl{}{}
D(G)^*\hfl{}{}\Theta\hfl{}{} 0$$
que l'on dualise
$$0\hfl{}{} D(G)\hfl{}{} D(E')\hfl{}{}
D(Q)\oplus\Theta\hfl{}{} 0\ ,$$
et H$^0(G)^*$ est un sous-espace vectoriel de H$^0(D(G))$
de dimension $n_G+l_G$.\par
\noindent Par construction, on a $n_{D(G)}=l_G$,
$n_{D(E')}=l$, $n_{D(Q)}=l_{Q}$ et:\par
- $d_{D(G)}=d_G=n_G+gl_G+l'_Q$\par
- $d_{D(Q)}= n_Q+gl_Q-l'_Q-d_{\Theta}$\medskip
\noindent De plus $D(G)$ est engendr\'e par ses sections et  
h$^1(D(G))=0 $ (cf lemme C-2 du Chap. 2). Maintenant soit
$F$ et $F'$ deux fibr\'es tels que $n_F=n_{D(G)} $,
$n_{F'}=n_{D(Q)} $, et
$d_F=d_{D(G)}$, $n_{F'}=n_{D(Q)}$, et $F$ est engendr\'e
par ses sections, et h$^1(F)=0$. Supposons que l'on ait une
suite exacte 
$$0\hfl{}{} F\hfl{}{} D(E')\hfl{}{}
F'\oplus\Theta\hfl{}{} 0\ .$$
Soit un sous-espace vectoriel $V$ de H$^0(F)$ de dimension 
$n_G+l_G$ qui engendre $F$. En posant 
$$G\simeq \bigl(ev:V\otimes{\cal O}\hfl{}{} F\bigr)^* $$
et 
$$Q\simeq {\rm Ker}\bigl(E'\hfl{}{} G\bigr)$$
on r\'eobtient une suite exacte
$$0\hfl{}{} Q\hfl{}{} E'\hfl{}{}
G\hfl{}{} 0\ .$$
Les suites exactes 
$$0\hfl{}{} Q\hfl{}{} E'\hfl{}{}
G\hfl{}{} 0$$
sont donc en correspondance bijective avec les suites
exactes
$$0\hfl{}{} F\hfl{}{} D(E')\hfl{}{}
F'\oplus\Theta\hfl{}{} 0$$
munies d'un sous-espace $V\subset{\rm H}^0(F)$ comme
ci-dessus. Pour $F$ fix\'e, le choix de $V$ est classifi\'e par
la grassmanienne Grass$_{n_G+l_G}\bigl({\rm H}^0(F)
\bigr)$ et donc d\'epend de $l'_Q(n_G+l_G)$ param\`etres 
(car h$^0(F)=n_G+l_G+l'_Q$). Il nous reste \`a compter le
nombre $\Gamma_1 $ de param\`etres  dont d\'ependent les
extensions
$$0\hfl{}{} F\hfl{}{} D(E')\hfl{}{}
F'\oplus\Theta\hfl{}{} 0$$
pour $D(E')$ fix\'e et on pourra poser
$$\Gamma=\Gamma_1+l'_Q(n_G+l_G)\ .$$
\noindent Pour calculer $\Gamma_1$, on suppose d'abord
que $D(E')$ varie dans $U_{l,d}$. Les familles des $F$ et
des $F'$ sont limit\'ees. Le choix du fibr\'e $F$
d\'epend de
$l_G^2(g-1)+1$ param\`etres; le choix du fibr\'e
$F'$ d\'epend de
$l_Q^2(g-1)+1$ param\`etres et le choix de $\Theta$
d\'epend de $n_{\Theta}$ param\`etres, le nombre de points
du support de ${\Theta}$. Pour
$F$ et
$F'$ fix\'es les extensions
$$0\hfl{}{} F\hfl{}{} D(E')\hfl{}{}
F'\oplus\Theta\hfl{}{} 0$$
sont param\'etr\'ees par
Ext$^1(F'\oplus \Theta, F) ={\rm H}^1(F'^*\otimes F)\oplus
{\rm Ext}^1(\Theta,F)$;
$D(E')$ devant \^etre stable, on a h$^0(F'^*\otimes
F)=0$. \par
\noindent Le
th\'eor\`eme de Riemann-Roch implique alors
que\medskip
\line{\hbox
to 3cm{\hfil${\rm h}^1(F'^*\otimes
F)$\hfil}$=l_G(n_Q+gl_Q-l'_Q-d_{\Theta})
-l_Q(n_G+gl_G+l'_Q+d_{\Theta})
+l_Ql_G(g-1)$\hfil}\par
\line{\hbox
to 3cm{}$=l_Gn_Q-l_Qn_G
+l_Ql_G(g-1)-l(d_{\Theta}+l'_Q)$\hfil}\par
\noindent De plus dim$\,\bigl({\rm Ext}^1(\Theta,F)\bigr)
=l_Gd_{\Theta}$. Donc, en consid\`erant
l'action de Aut$(\Theta)$ sur ces extensions (dim$\bigl(
{\rm Aut}(\Theta)\bigr)\geq n_{\Theta}$), on en d\'eduit que
le nombre de param\`etres dont d\'ependent les extensions
ci-dessus pour $D(E')$, $F$, $F'$ et $\Theta$ variants (mais
les rangs et les degr\'es sont fix\'es) est donc inf\'erieur
ou \'egal \`a
$$ l_G^2(g-1)+1+l_Q^2(g-1)+1+l_Gn_Q-l_Qn_G
+l_Ql_G(g-1)-(l-l_G)d_{\Theta}-ll'_Q-1\ ;$$
apr\`es simplification et avec $l=l_G+l_Q$, l'expression
ci-dessus est encore \'egale \`a
$$ l^2(g-1)+1+l_Gn_Q-l_Qn_G
-l_Ql_G(g-1)-l_Qd_{\Theta}-ll'_Q\ .$$
\noindent Or
$D(E')$ d\'epend de
$l^2(g-1)+1$ param\`etres, donc pour  $D(E')$ fix\'e
assez g\'en\'eral, les suites exactes comme
ci-dessus doivent d\'ependre de
 au plus
$$l_Gn_Q-l_Qn_G
-l_Ql_G(g-1)-l_Qd_{\Theta}-ll'_Q$$
param\`etres et l'on posera en fait
$$\Gamma_1=l_Gn_Q-l_Qn_G
-l_Ql_G(g-1)-ll'_Q\ .$$
L'in\'egalit\'e $(*)$ que l'on veut d\'emontrer s'\'ecrit
$$l'(n_G+l_G)>l_Gn_Q-l_Qn_G
-l_Ql_G(g-1)-ll'_Q+l'_Q(n_G+l_G)$$
ou encore
$$(l'-l'_Q)n_G>l_Gn_Q-l_Qn_G
-l_Ql_G(g-1)-l'_Ql_Q-l'l_G\ .$$
\noindent Or, $\mu(H)={d_Q\over n_Q-l'}=
{n_Q+gl_Q-l'_Q\over n_Q-l'}\geq 1+{gl_G+l'_Q\over
n_G}=\mu(G)$ implique que
$$n_G(l'+gl_Q-l'_Q)\geq n_Q(gl_G+l'_Q)-l'(gl_G+l'_Q)$$
 donc que
$$l_Gn_Q-l_Qn_G\leq {1\over g}\lbrack
(l'-l'_Q)n_G-l'_Qn_Q\rbrack+l'l_G+{l'l'_Q\over g}\ ,$$
et il nous suffira de prouver
$${g-1\over g}\lbrack
(l'-l'_Q)n_G\rbrack>-{l'_Q\over g}n_Q+ l'l_G+{l'l'_Q\over
g}-l_Ql_G(g-1)-l'_Ql_Q-l'l_G\ .$$
Mais on a $n_G> gl_G$ (cf Th\'eor\`eme A-1 du
Chap. 2), $l'_Q<l'<g $ et $l_G>0 $ (car $G$ engendr\'e par
ses sections), donc
$${g-1\over g}\lbrack
(l'-l'_Q)n_G\rbrack>l'l_G $$
ainsi que
$${l'l'_Q\over
g}-l'l_G<0 $$
ce qui montre que l'in\'egalit\'e  $(*)$ est toujours vraie.\par 
\noindent On vient donc de montrer que si le fibr\'e $Q$ est
g\'en\'eriquement engendr\'e par ses sections pour
"$Q$ assez g\'en\'eral", alors g\'en\'eriquement un sous-espace
de sections
${\cal O}^{l'}\in {\rm Grass}_{l'}\bigl({\rm
H}^0(E')\bigr)$  ne s'envoie dans aucun fibr\'e
$Q$ comme ci-dessus et donc le conoyau $E$ de l'injection 
serait stable.\par
\noindent  Pour  $l'=1$, on a fini.\bigskip
\noindent {\bf La r\'ecurrence:} L'hypoth\`ese de
r\'ecurrence est: 
Soit $E'$ un fibr\'e g\'en\'erique comme
ci-dessus et soit
${\cal O}^{l'-1}\subset {\rm H}^0(E')$ un
sous-espace vectoriel de dimension
$l'-1$ g\'en\'erique dans la grassmanienne,
Grass$_{l'-1}\bigl({\rm H}^0(E')\bigr)$, des
sous-espaces vectoriels de dimension
$l'-1$ dans
${\rm H}^0(E')$. Alors, le conoyau $E$
$$0\hfl{}{} {\cal O}^{l'-1}\hfl{}{} E'\hfl{}{} E\hfl{}{}0$$
est stable.\par
\noindent Maintenant, soit ${\cal O}^{l'}
\subset {\rm H}^0(E')$ un
sous-espace vectoriel de dimension
$l'$ g\'en\'erique dans la grassmanienne,
Grass$_{l'}\bigl({\rm H}^0(E')\bigr)$,
des
sous-espaces vectoriels de dimension $l'$ dans
${\rm H}^0(E')$. On veut montrer que le
conoyau $E$
$$0\hfl{}{} {\cal O}^{l'}\hfl{}{} E'\hfl{}{} E\hfl{}{}0$$
est stable.\par
\noindent On a vu que la non-stabilit\'e de $E$
impliquait l'existence d'une suite exacte
$$0\hfl{}{} H\hfl{}{} E\hfl{}{} G\hfl{}{}0$$
o\`u $G$ est stable de pente $\leq \mu(E)$ et 
donne, par le Diagramme 1, une suite exacte
$$0\hfl{}{} Q\hfl{}{} E'\hfl{}{} G\hfl{}{}0$$
D'apr\`es ce qui pr\'ec\`ede, il nous reste \`a montrer que si $E'$
et ${\cal O}^{l'}$ sont assez g\'en\'eraux, alors $Q$ est 
g\'en\'eriquement engendr\'e par ses sections\par

\noindent On peut supposer que g\'en\'eriquement
un sous-espace  de sections ${\cal O}^{l'-1}
\subset {\cal O}^{l'}$ v\'erifie: le conoyau $R$
$$0\hfl{}{} {\cal O}^{l'-1}\hfl{}{} E'\hfl{}{} R\hfl{}{}0$$
est stable.\par
\noindent Cette suite exacte se place dans un
diagramme commutatif
$$\matrix{
&&&&&&0&&&\cr
&&&&&&\downarrow&&&\cr
&&0&&0&\longrightarrow&{\cal O}&&&\cr
&&\downarrow&&\downarrow&&\downarrow&&&\cr
0&\longrightarrow&{\cal O}^{l'-1}&\longrightarrow& 
E'&\longrightarrow
&R &\longrightarrow&0\cr
&&\downarrow&&\downarrow&&\downarrow&&&\cr
0&\longrightarrow&{\cal O}^{l'}&\longrightarrow&
E'&
\longrightarrow& E&\longrightarrow&0\cr
&&\downarrow&&\downarrow&&\downarrow&&&\cr
&&{\cal O}&\longrightarrow&0
&
&0&&\cr
&&\downarrow&&&&&&&\cr
&&0&&&&&&&\cr
}$$
dont on retient la suite exacte
$$0\hfl{}{} {\cal O}\hfl{}{} R\hfl{}{} E\hfl{}{}0$$
qui, avec la suite exacte en $E$, $G$ et $H$,
donne un nouveau diagramme commutatif
$$\matrix{
&&&&&&0&&&\cr
&&&&&&\downarrow&&&\cr
&&0&& 
0&\longrightarrow
&H &&\cr
&&\downarrow&&\downarrow&&\downarrow&&&\cr
0&\longrightarrow&{\cal O}&\longrightarrow&
R&
\longrightarrow& E&\longrightarrow&0\cr
&&\downarrow&&\downarrow&&\downarrow&&&\cr
0&\longrightarrow&Q_0&\longrightarrow&R
&\longrightarrow
&G&\longrightarrow&0\cr
&&\downarrow&&\downarrow&&\downarrow&&&\cr
&&H&\longrightarrow&0&&0&&&\cr
&&\downarrow&&&&&&&\cr
&&0&&&&&&&\cr}$$
\centerline{\it Diagramme 3}\medskip
\noindent o\`u $Q_0$ est le noyau du morphisme
$R\rightarrow G$\par
\noindent Comme $\mu(H)\geq\mu(E)>\mu(R)$,
on en d\'eduit que h$^0(Q_0^*)=0$ (sinon, $H$
serait un sous-fibr\'e de $R$ stable). On a encore
un dernier diagramme
$$\matrix{
&&&&&&0&&&\cr
&&&&&&\downarrow&&&\cr
&&0&& 
0&\longrightarrow
&Q_0 &&\cr
&&\downarrow&&\downarrow&&\downarrow&&&\cr
0&\longrightarrow&{\cal O}^{l'-1}&\longrightarrow&
E'&
\longrightarrow& R&\longrightarrow&0\cr
&&\downarrow&&\downarrow&&\downarrow&&&\cr
0&\longrightarrow&Q&\longrightarrow&E'
&\longrightarrow
&G&\longrightarrow&0\cr
&&\downarrow&&\downarrow&&\downarrow&&&\cr
&&Q_0&\longrightarrow&0&&0&&&\cr
&&\downarrow&&&&&&&\cr
&&0&&&&&&&\cr
}$$
dont on retient la suite exacte
$$0\hfl{}{} {\cal O}^{l'-1}\hfl{}{} Q\hfl{}{}
Q_0\hfl{}{}0\ .$$
On en d\'eduit que
$\mu(Q_0)=1+{gl_{Q}+l'-l'_Q-1\over n_{Q_0}}$. Comme
h$^0(Q_0^*)=0$ et que $Q_0$ est un sous-fibr\'e de $R$
stable de pente $<2$, la proposition A-2 du chap. 1 implique
que
 h$^0(Q_0)=n_{Q_0}+l_Q$. Toutes les sections de $Q_0$
remontent donc \`a $Q$. On va montrer que $Q_0$ est
g\'en\'eriquement engendr\'e par ses sections et ce sera
fini.\par
\noindent
C'est  la suite exacte donn\'ee
par le diagramme 3:
$$0\hfl{}{} {\cal O}\hfl{}{} Q_0\hfl{}{}
H\hfl{}{}0$$
qui nous permet de conclure:
d'apr\`es ce qui pr\'ec\`ede, on a
$\mu(H)=1+{gl_{Q}+l'-l'_Q\over n_{H}}$. 
et comme $H$ contredit la stabilit\'e de $E$,
on a aussi par hypoth\`ese
$\mu(H)>\mu(R)$. Comme dans le cas $l'=1$, on
en d\'eduit que $Q_0$ est g\'en\'eriquement
engendr\'e par ses sections (sinon $Q_0$ contient
un sous-faisceau de pente $\geq 1+{g(l_{Q}+1)\over
n_{H}}$ et ceci contredit la stabilit\'e de $R$).\medskip
\noindent La d\'emonstration du
th\'eor\`eme est donc termin\'ee.\hbox to
1cm{}$\diamondsuit$\vfill\break{\nopagenumbers\line{}\vskip 2cm
\line{\rmgXV Annexe:\hfill}\vskip 2cm
\centerline{\rmgXV Les figures}

\vfill\break}\line{}\par 

\vfill\break\line{}\par 

\vfill\break\line{}\par 

\vfill\break\line{}\par 

\vfill\break\line{}\par 

\vfill\break\line{}\par 

\vfill\break

\centerline{\rmgXIII Bibliographie}\vskip 2cm

\cro{A}{Atiyah, M.F.}{Vector bundles on an elliptic curve.}
Proc. Lond. Math.Soc. (3) {\bf 7}. (1982) pp. 414-452.\medskip

\cro{A-C-G-H}{Arbarello, E., Cornalba, M., Griffiths,
P., Harris, J.}{Geometry of Algebraic
curves I.} Grundlehren Math. Wiss. {\bf 267}, Springer verlag, New-York, (1984).\medskip

\cro{B-G-N}{Brambila-Paz, L., Grzegorczyk, I., Newstead, P.
E.}{Geography of Brill-Noether loci for small slopes} Preprint
Octobre 95.\medskip

\cro{B-N}{Brambila-Paz, L, Newstead, P. E.}{Subvariedades
del espacio moduli} Memorias del XXVII Congreso de la
Sociedad Matematica Mexicana, Comunicaciones,
Aportaciones Matematicas {\bf 16} (1995), pp. 43-53.
\medskip

\cro{B}{Butler C., D.}{Normal generation of vector bundles
over a curve.} J. Differential Geometry {\bf 39}
(1994) pp. 1-34.\medskip

\cro{EGA}{Grothendieck, A.}{El\'ements de G\'eom\'etrie
Alg\'ebrique} Springer Verlag {\bf 166} (1971)
\medskip

\cro{E-H-1}{Eisenbud, D., Harris, J.}{Limit linear series:
Basic theory.}
Invent. Math. {\bf 85} (1986) pp. 337-371.\medskip

\cro{E-H-2}{$\vcenter{\hrule width 3cm}$}{Divisors on
general curves and cupsidal rational
curves.}  Invent. Math. {\bf 74} (1983) pp. 371-418.\medskip

\cro{E-H-3}{$\vcenter{\hrule width 3cm}$}{A simpler proof of
Gieseker-Petri Theorem on special divisors.} Inventiones
math. {\bf 74} (1983) pp. 269-280.\medskip

\cro{G}{Grothendieck, A.}{Technique de descente et
th\'eor\`eme d'existence en g\'eom\'etrie
alg\'ebrique. IV Les sh\'emas de Hilbert.}
S\'eminaire Bourbaki {\bf 221}, (1960/61).\medskip

\cro{H}{Hartshorne, R.}{Algebraic geometry.} 
Berlin-Heidelberg-New York Springer Verlag (1977) . \medskip

\cro{L}{Laumon, G.}{Fibr\'es vectoriels sp\'eciaux.} Bull. Soc.
Math. France {\bf 119} (1991) pp. 97-119.\medskip

\cro{Li}{Li, Y.}{Spectral curves, theta divisors and Picard
bundles}, International J. Math. 2 (1991), pp.
525-550.\medskip

\cro{LP}{Le Potier, J.}{Fibr\'es vectoriels sur les
courbes alg\'ebriques.},
Publications Math\'ematiques de l'universit\'e Paris 7- Denis
Diderot {\bf 35} (1996) .\medskip 

\cro{M-F}{Mumford, D., Fogarty, J.}{Geometric Invariant
Theory.}
Springer Verlag (1982) .\medskip

\cro{Na-R}{Narasimhan, M. S., Ramanan, S.}{Deformation of
the moduli space of vector bundles over an algebraic curve}
Ann.-Math. (2) {\bf 101} (1975), pp. 391-417.\medskip 

\cro{Na-Se}{Narasimhan, M. S., Seshadri, C. S.}{Stable and
unitary vector bundles on a compact Riemann surface}
Ann.-of-Math. (2) {\bf 82} (1965) pp. 540-567.\medskip

\cro{Pe}{Peskine, C.}{Introduction alg\'ebrique \`a la
g\'eom\'etrie projective} Cours de D.E.A (1993).\medskip

\cro{P-R}{Paranjape, P., Ramanan, S.}{On the canonical ring
of a curve} Collection: Algebraic geometry and commutativ
algebra, Vol. II, (1988) pp. 503-516\medskip

\cro{Re}{Re, R.}{Multiplication of sections and Clifford
bounds for special stable bundles on curves} Preprint
96.\medskip

\cro{Rego}{Rego, C. J.}{Deformation of modules on curves
and surfaces.} Singularities, Representatation of Algebras,
and Vector Bundles. Proceedings, Lambrecht 1985, Springer
Lect. Notes in Math. {\bf 1273} pp. 157-167\medskip

\cro{Se}{Seshadri, C. S.}{Fibr\'es vectoriels sur les
courbes alg\'ebriques.}
Ast\'erisque {\bf 96}. (1982) .\medskip

\cro{Su}{Sundaram, N.}{Special divisor and vector
bundles} T\^ohoku Math. Journ. {\bf 39} (1987)
pp.175-273.\medskip

\cro{T}{Tan, X. J.}{Some results on the existence of rank two
special stable vector bundles} Manuscripta Math. {\bf 75}
(1992), 365-373.\medskip

\cro{Te-1}{Monserrat Teixidor I Bigas}{Brill-Noether
Theory for stable vector bundles.},
Duke math. Journal {\bf 62} (1991) .\medskip

\cro{Te-2}{$\vcenter{\hrule width 3cm}$}{Brill-Noether
Theory for vector bundles of rank 2.} 
T\^ohoku Math. journal {\bf 43} (1991) pp. 123-126.\medskip

\cro{Te-3}{$\vcenter{\hrule width 3cm}$}{On the
Gieseker-Petri map for rank $2$ vector bundles}
Manuscripta Math. {\bf 75} (1992) pp. 375-382.\medskip
\vfill\break

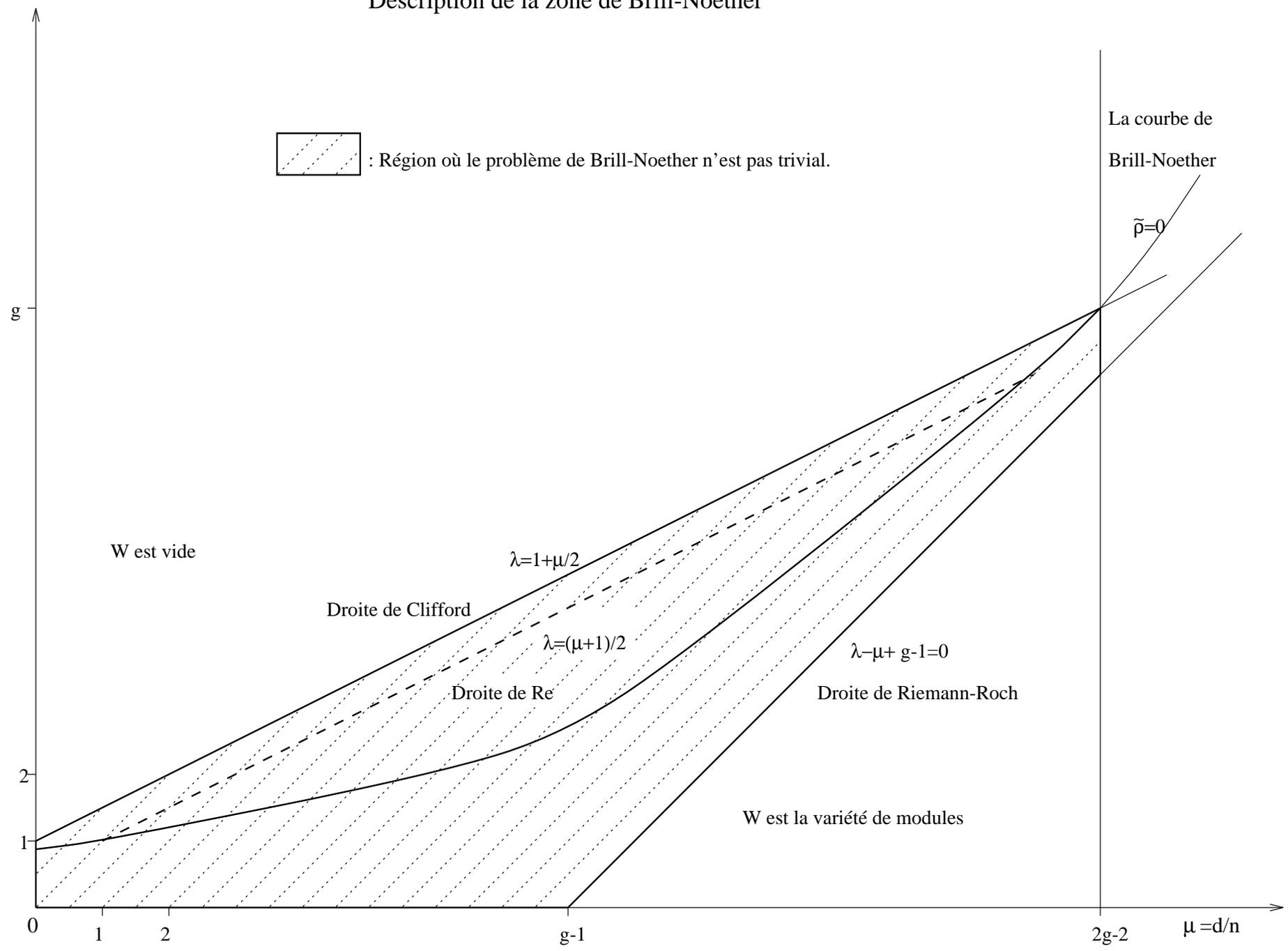
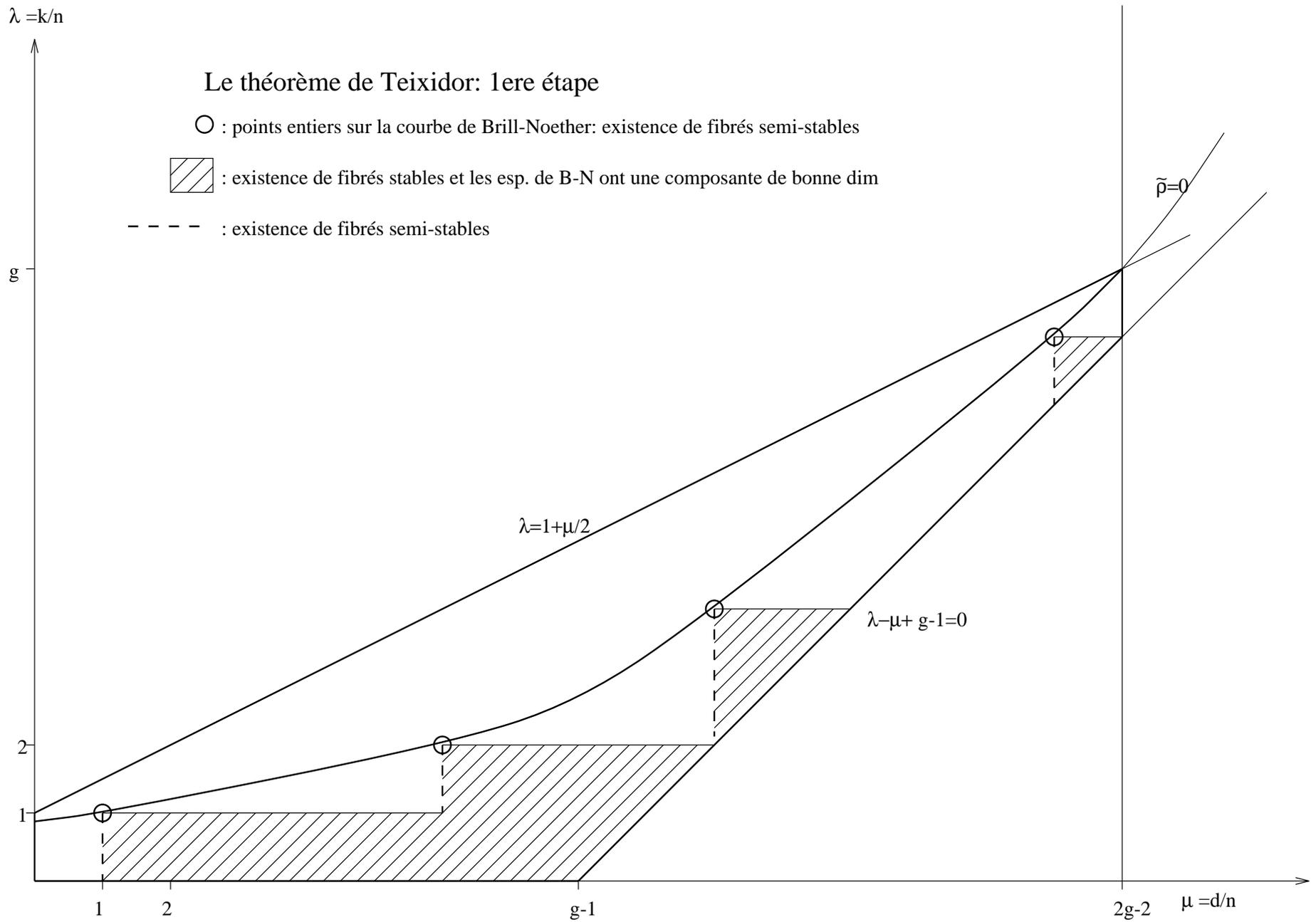
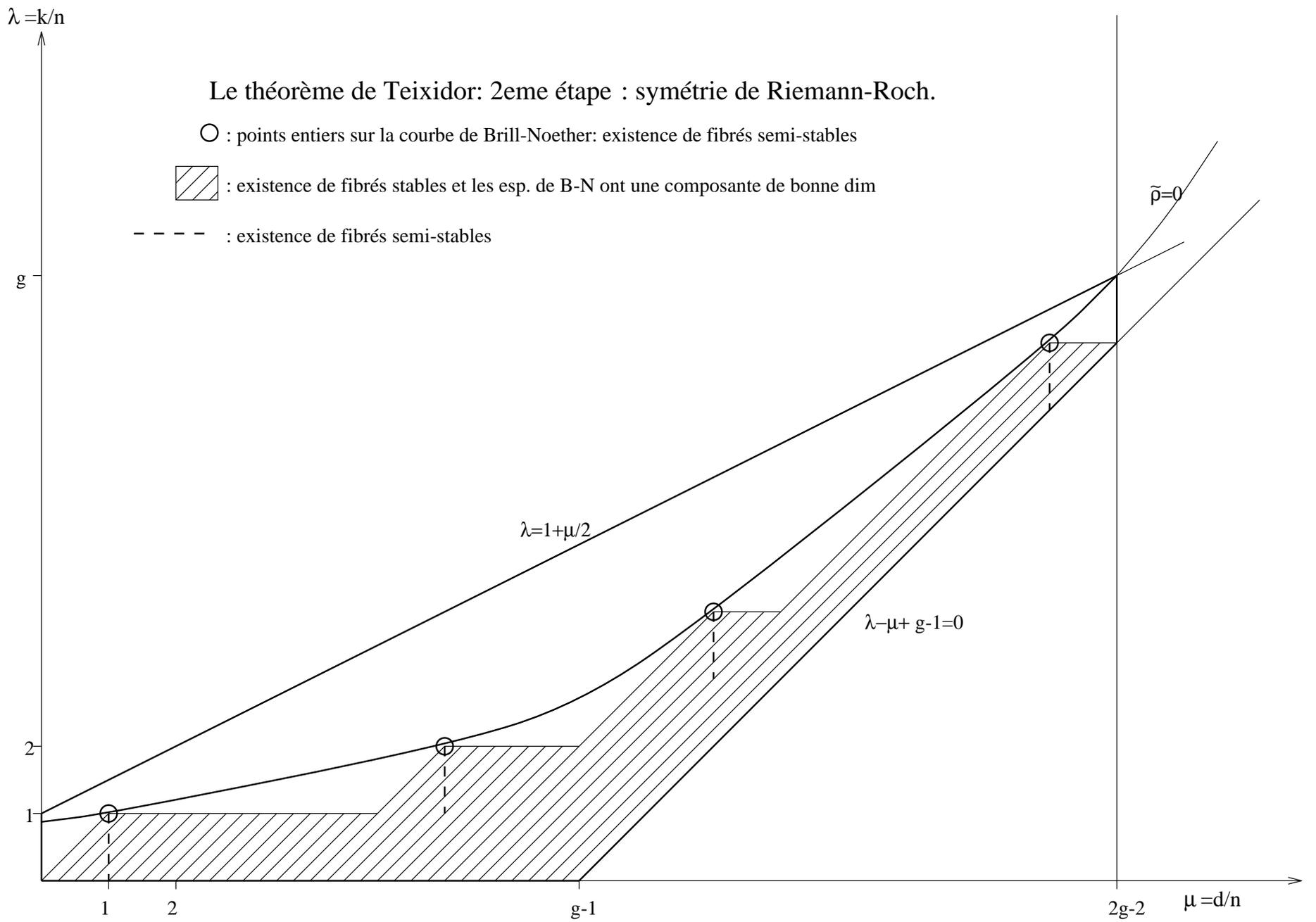
\end